\documentclass[a4paper,10pt]{amsart}
\usepackage{epsf}  
\usepackage{amsfonts, amssymb, amsthm, amsmath}  
\usepackage{hyperref}
\usepackage{tikz}
\usetikzlibrary{matrix,arrows,decorations.pathmorphing}
\usepackage{tikz-cd}
\usepackage{graphicx}
\DeclareMathOperator{\ev}{\mathfrak{ev}}
\usepackage{psfrag}
\newtheorem{thm}{Theorem}

\newtheorem{lemma}[thm]{Lemma}  
\newtheorem{remark}[thm]{Remark}  
\newtheorem{defn}[thm]{Definition}  
\newtheorem{prop}[thm]{Proposition}  
\newtheorem{claim}[thm]{Claim}  
  
\numberwithin{thm}{section}  
\def\pf{\noindent\emph{Proof: }}  
\def\stop{\hfill$\square$}  
\def\co{\colon\thinspace}

\providecommand{\totl}[1]{\ensuremath{\lceil #1\rceil }}
\providecommand{\totb}[1]{\ensuremath{\underline{ #1}}}
\DeclareMathOperator{\End}{\mathcal E}
\DeclareMathOperator{\gend}{\tilde{\mathcal E}}
\DeclareMathOperator{\rend}{End}
\DeclareMathOperator{\Aut}{Aut}
\DeclareMathOperator{\cut}{cut}

\newcommand{\rh}{{}^{r}H}
\newcommand{\ie}{\text{ied}}

\newcommand{\ex}{\bold}
\providecommand {\e}[1]{\mathfrak t^{#1}}

\providecommand{\C}[2]{\ensuremath {C^{#1,\underline{#2}}}}

\newcommand{\tc}[1]{\check\rvert_{#1}}

\newcommand{\Msw}{\mathcal M^{st}}
\newcommand{\Mpsw}{\mathcal M^{\infty,\underline 1}}
\newcommand{\Mu}{\mathcal M^{\bullet}}

\DeclareMathOperator{\id}{id}
\DeclareMathOperator{\expl}{Expl}

\newcommand{\dbar}{\bar{\partial}}

\providecommand{\et}[2]{\ensuremath{\bold T^{#1}_{#2}}}
\providecommand{\lrb}[1]{\ensuremath{\left(#1\right)}}
\providecommand{\abs}[1]{\left\lvert #1\right\rvert}

\author{Brett Parker   }
\email{brettdparker@gmail.com}  
\thanks{Partially supported by ARC grant  DP140100296.}
  
\title{Tropical gluing formulae for Gromov--Witten invariants}

\begin{document}
\maketitle

\begin{abstract}We prove two tropical gluing formulae for Gromov--Witten invariants of exploded manifolds, useful for calculating Gromov--Witten invariants of a symplectic manifold using a normal-crossing  degeneration. The first formula  generalizes   the symplectic-sum formula for Gromov--Witten invariants. The second formula is stronger, and also generalizes Kontsevich and Manin's splitting and genus-reduction axioms. Both tropical gluing formulae  have versions incorporating gravitational descendants.
\end{abstract}

\tableofcontents

\newpage

\section{Introduction}

This paper explains the following  tropical gluing formulae for Gromov--Witten invariants. 
\begin{equation}\label{gf1}\eta \tc \gamma=\frac{k_{\gamma}}{\abs{\Aut\gamma}}i^{[\gamma]}_{!}\Delta^{*}\prod_{v}\eta^{[\gamma_{v}]}\end{equation}
\begin{equation}\label{gf2}\mu\tc \gamma=\frac 1{\abs{\Aut\gamma}}I_{!}^{[\gamma]}\Delta^{*}\prod_{v}\mu^{[\gamma_{v}]} \end{equation}
On the left,  $\eta$ and $\mu$ represent Gromov--Witten invariants. In the case of of a compact symplectic manifold $B$, $\eta$ and $\mu$ encode the usual Gromov--Witten invariants obtained using evaluation of curves with $n$ marked points to $B^{n}$ and $B^{n}\times \bar M_{g,n}$ respectively.  The notation $\tc \gamma$ indicates the contribution of a tropical curve $\gamma$ to these invariants. On the right, $\eta^{[\gamma_{v}]}$ and $\mu^{[\gamma_{v}]}$ indicate relative\footnote{These relative invariants are defined using exploded manifolds, however they are roughly equivalent to Ionel's relative invariants from \cite{IonelGW}, and in the algebraic setting are equivalent to log Gromov--Witten invariants, \cite{GSlogGW,Chen,acgw}.} Gromov--Witten invariants associated to vertices $v$ of $\gamma$, and the remaining terms are elementary instructions for combining these relative invariants to compute $\eta$ and $\mu$. The relative invariants themselves are  subject to the same gluing formulae, which often suffice to determine them; see for example \cite{tec,3d}.

Formula (\ref{gf1}) generalizes the symplectic-sum formula for Gromov--Witten invariants, \cite{Li,IP,ruan}, to normal-crossing or log-smooth degenerations, and symplectic analogues;\footnote{For an appropriate symplectic analogue of normal-crossing divisors, see \cite{ZingerNC}, and for a symplectic analogue of a normal-crossing degeneration, see \cite{symplecticMS}.} see \cite{tpgf} for an exposition of the simplest generalization ---  the case of a triple-product. Formula (\ref{gf2}) also gives a degeneration formula for Gromov--Witten invariants in these situations, but includes the contribution of cohomology classes from Deligne-Mumford space; without applying any degeneration, it 
 specializes to Kontsevich and Mannin's splitting and genus-reduction axioms for Gromov--Witten invariants of symplectic manifolds. We also prove a version of each formula incorporating  gravitational descendants, the contribution of Chern classes of tautological line-bundles over the moduli stack of curves.

\

This paper's action takes place within the category of exploded manifolds;\footnote{See \cite{iec} or \cite{scgp} for an introduction to exploded manifolds. In the algebraic setting, using exploded manifolds is almost equivalent to using log schemes,  and equivalent gluing formulae should be provable entirely within the framework of log schemes. For construction of log Gromov--Witten invariants, see \cite{GSlogGW,Chen,acgw}. See \cite{elc} for a log-exploded  dictionary.} our results apply to compact symplectic manifolds using degeneration. We can explode any normal-crossing degeneration of a symplectic manifold $B$ to produce a smooth family of exploded manifolds;  smooth fibers of the original degeneration are unchanged, but the singular fiber is replaced with a family of exploded manifolds. Gromov--Witten invariants do not change in families of exploded manifolds, so the Gromov--Witten invariants of $B$ can be computed using any other fiber, and it is convenient to use one of the fibers $\ex B$  replacing the singular fiber. For examples, see \cite{3d,tpgf,tec}. 

Let us explain the tropical curves $\gamma$ appearing in our gluing formulae, (\ref{gf1}) and (\ref{gf2}). Each exploded manifold $\ex B$ has a tropical part, $\totb{\ex B}$, consisting of a complex of integral-affine polytopes. For example, if $\ex B$ arises as a fiber of a normal-crossing degeneration with singular fiber $B'$,  $\totb{\ex B}$ will be the dual complex of simplices with a vertex for each component of $B'$, and a $n$--simplex for every intersection of $(n+1)$ components. Each curve in $\ex B$ is a map $f\co  \ex C\longrightarrow \ex B$ which itself has a tropical part, $\totb f\co  \totb{\ex C}\longrightarrow \totb{\ex B}$, where $\totb{\ex C}$ is a graph with a complete integral-affine structure on its edges, and $\totb f$ is integral-affine. A tropical curve $\gamma$ in $\totb{\ex B}$ is such an integral-affine map from a complete integral-affine graph.

 All exploded manifolds $\ex B$ also have a smooth part, $\totl{\ex B}$ consisting of a complex of smooth manifolds. In the above case, $\totl{\ex B}$ is isomorphic to the singular fiber $B'$, and $\totl f\co  \totl{\ex C}\longrightarrow \totl{\ex B}$ is a holomorphic curve in $B'$, with a component for each vertex of $\totb{\ex C}$, a node for each internal edge of $\totb{\ex C}$, and a marked point corresponding to each end of $\totb{\ex C}$ (isomorphic to $[0,\infty)$). Although the reader might  intuitively think in terms of these smooth parts, consideration of the extra tropical structure is essential for correct gluing and degeneration formulae. 

 Gromov--Witten invariants of $\ex B$ decompose into a sum of (virtual) contributions of tropical curves $\gamma$, or rather the  holomorphic curves with tropical part isomorphic to $\gamma$. Our gluing formulae compute this contribution of $\gamma$ in terms of relative invariants, $\eta^{[\gamma_{v}]}$ and $\mu^{[\gamma_{v}]}$. These relative invariants are Gromov--Witten invariants of an exploded manifold $\ex B\tc v$ created by completing the stratum of $\ex B$ containing $v$  using the tropical completion described in section 7 of \cite{vfc}. In the case that $\totl{\ex B}$ is the degenerate fiber of a normal-crossing degeneration and $v$ is in a corner of $\totb{\ex B}$, the corresponding stratum of $\totl{\ex B}$ is a manifold with a normal-crossing divisor, and $\ex B\tc v$ is the explosion\footnote{See \cite{iec},\cite{scgp}, or \cite{elc} for an exposition of the explosion functor.} of this stratum. Our relative invariants $\eta^{[\gamma_{v}]}$ and $\mu^{[\gamma_{v}]}$ are Gromov--Witten invariants of this stratum relative to its normal-crossing divisor. If $v$ is in a $k$--dimensional stratum of $\totb{\ex B}$, $\ex B\tc v$ is a   $k$--complex-dimensional bundle over the corresponding $k$--fold intersection of components of $\totl{\ex B}$, related to the `rubber components' or `expansions' that appear in \cite{Li,IP,ruan}. The Gromov--Witten invariants of such `expansions' are not important for the symplectic-sum formula, but are critical for its correct generalization; for simple examples  in the case of a triple-product, see \cite{tpgf}.
 
 Let us describe our  Gromov--Witten invariant $\eta$.
 \[\eta:=ev_{!}(q^{E}\hbar^{2g-2+n})\]
 Above,  $q$ and $\hbar$ are dummy variables whose exponents record the symplectic energy $E$, and Euler characteristic, $2g-2+n$, of curves, and $ev_{!}$ indicates pushforward\footnote{See \cite{vfc}, section 5.3,  for pushforwards from the virtual fundamental class.} from the virtual fundamental class\footnote{See \cite{evc} for the  construction of an embedded Kuranishi structure for the moduli stack of holomorphic curves, and  \cite{vfc} for the construction of the virtual fundamental class using this embedded Kuranishi structure.} using the  evaluation map,
 \[ev\co  \Msw_{\cdot}(\ex B)\longrightarrow \coprod_{n}(\rend \ex B)^{n}\]
constructed in section \ref{simple ev section}. The notation $\Msw_{\cdot}(\ex B)$ indicates a moduli stack of (not-necessarily holomorphic) stable curves in $\ex B$ with labeled ends\footnote{An end of an (exploded) holomorphic curve $\ex C$  corresponds to an end, or infinite edge, of its underlying tropical curve $\totb{\ex C}$. Ends correspond to marked points of the underlying nodal curve $\totl{\ex C}$.}.
When $\totb{\ex B}$ is bounded, $\rend\ex B=\ex B$, and our evaluation map $ev$ simply records the location of ends/punctures of curves. When $\ex B$ is the explosion of a manifold with a  smooth divisor $D$, $\rend \ex B$ is the disjoint union of $\ex B$ with a copy of $D$ for each positive integer `contact order'.\footnote{The tropical part of $\ex B$ in this case is a half-line, as is the tropical part of each end of a curve. The smooth part of a curve in $\ex B$  has contact order $k$ with $D$ at an end if the derivative of its tropical part at that end is $k$.} The evaluation map at an end/puncture  lands in $\ex B$ or the appropriate copy of $D$, depending on the contact order at that end. In more general cases,\footnote{The evaluation map $ev$ is described precisely in section \ref{simple ev section}.
}  $\rend \ex B$ records the possible positions of the ends of curves in $\ex B$, and is the disjoint union of $\ex B$ with exploded manifolds of real-dimension $2$ lower than $\ex B$.  

The relative invariant $\eta ^{[\gamma_{v}]}$ appearing in our gluing formula (\ref{gf1}) is the restriction of the corresponding invariant $\eta$ of $\ex B\tc v$ to the connected component
\[\rend^{[\gamma_{v}]}\ex B\tc v\subset \coprod_{n}(\rend \ex B\tc v)^{n}\]
containing the image of curves in $\ex B\tc v$ with tropical part $\gamma_{v}$, where $\gamma_{v}$ is a tropical curve in $\totb{\ex B\tc v}$ with a single vertex, and edges created by infinitely lengthening all edges of $\gamma$ leaving $v$.  Alternatively, we may define
\[\eta^{[\gamma_{v}]}:=ev^{[\gamma_{v}]}_{!}(q^{E}\hbar^{2g-2+n})\]
where 
\[ ev^{[\gamma_{v}]}\co  \Msw_{[\gamma_{v}]}(\ex B\tc v)\longrightarrow \rend^{[\gamma_{v}]}(\ex B\tc v)\]
is our evaluation map $ev$ restricted to the moduli stack $\Msw_{[\gamma_{v}]}(\ex B\tc v)\subset \Msw_{\cdot}(\ex B\tc v)$ of curves in $\ex B\tc v$ with ends labelled by the ends of $\gamma_{v}$, and having derivatives\footnote{If $f$ is a curve in $\ex B\tc v$, and $\ex B\tc v$ is the explosion of a manifold with normal-crossing divisor $D$, the derivative of an end of $\totb f$ is an integral-vector whose components correspond to contact orders of $\totl f$ with components of $D$.} at these ends equal to the derivatives of the ends of $\gamma_{v}$.  

Suppose that $\gamma$ has $n$ ends. To write our gluing formula, the relationship between $(\rend \ex B)^{n}$ and $\rend ^{[\gamma_{v}]}\ex B\tc v$ is encapsulated in the following maps, explained in section \ref{simple ev section}.
\[\begin{tikzcd}\prod_{v}\rend^{[\gamma_{v}]}(\ex B\tc v)&\lar{\Delta} \ex Y \rar{i^{[\gamma]}}&(\rend\ex B)^{n}\tc{\totb{ev}\gamma}\end{tikzcd}\]
Each internal edge of $\gamma$ corresponds to two ends of $\coprod_{v}\gamma_{v}$. The corresponding two connected components of $\coprod_{v}\rend (\ex B\tc v)$ are naturally isomorphic, and the map $\Delta$ is the inclusion of the  diagonal subset of $\prod_{v}\rend^{[\gamma_{v}]}(\ex B\tc v)$ using these isomorphism for each internal edge of $\gamma$. The map $i^{[\gamma]}$ is a projection which forgets the information from each of these internal edges. In particular, each of the $n$ ends of $\gamma$ corresponds to a unique end of $\coprod_{v}\gamma_{v}$, and the connected component of $(\coprod_{v}\rend(\ex B\tc v))^{n}$ recording the position of these $n$ ends is naturally isomorphic to $(\rend\ex B)^{n}\tc{\totb{ev}\gamma}$, the tropical completion of $(\rend \ex B)^{n}$ at the image of curves with tropical part $\gamma$. Our map $i^{[\gamma]}$ is the projection which forgets the factors recording the position of internal edges, followed by this natural isomorphism. 

Our first gluing formula
\[\eta \tc \gamma=\frac{k_{\gamma}}{\abs{\Aut\gamma}}i^{[\gamma]}_{!}\Delta^{*}\prod_{v}\eta^{[\gamma_{v}]}\]
contains two combinatorial factors we have yet to explain. The constant $k_{\gamma}$ is the product of the multiplicities $m_{e}$ of the internal edges of $\gamma$, where each internal edge of $\gamma$ has derivative an integral-vector equal to $m_{e}$ times a primitive integral-vector. This factor of $k_{\gamma}$ arrises because our gluing formula follows from a natural fiber-product diagram involving not the exploded manifold $\rend \ex B$, but a corresponding stack $\End\ex B$ that is the quotient of $\rend\ex B$ by a trivial group action, $\mathbb Z_{m_{e}}$ on the component corresponding to an edge of multiplicity $m_{e}>0$, and the infinite group $\ex T$ on the component corresponding to an edge of multiplicity $0$. Our formula also requires division by the size of the automorphism group,\footnote{We only use automorphisms of $\gamma$ that fix ends, because we have labeled ends of curves by working with $\Msw_{\cdot}$.} $\Aut \gamma$, of the tropical curve $\gamma$, because the natural fiber-product diagram we use involves labelling edges of curves by the edges of $\gamma$. 

We also prove a related gluing formula that includes gravitational descendants. Suppose that $W$ is a tautological vectorbundle over $\Msw_{\cdot}(\ex B)$, so $W$ is a product of tautological line-bundles corresponding to the ends of curves. Let $W_{v}$ be the corresponding tautological vectorbundle over $\Msw_{[\gamma_{v}]}(\ex B\tc v)$, remembering that some ends of curves in $\Msw_{[\gamma_{v}]}$ correspond naturally to ends of curves in $\Msw_{\cdot}$. Define
\[\eta(W):=ev_{!}(q^{E}\hbar^{2g-2+n}c(W))\]
\[\eta^{[\gamma_{v}]}(W_{v}):=ev^{[\gamma_{v}]}_{!}(q^{E}\hbar^{2g-2+n}c(W_{v}))\]
using the pushforward of top Chern-classes $c(W)$ and $c(W_{v})$ defined in \cite{vfc}, Remark 5.2. These Gromov--Witten invariants satisfy the following modified version of (\ref{gf1}).
\[\eta(W) \tc \gamma=\frac{k_{\gamma}}{\abs{\Aut\gamma}}i^{[\gamma]}_{!}\Delta^{*}\prod_{v}\eta^{[\gamma_{v}]}(W_{v})\]

Let us consider the elements of our gluing formula in a quick example, discussed in section 8 of \cite{tpgf}. Degenerate $\mathbb CP^{2}$ into three components $M_{i}$ isomorphic to  $CP^{2}$ blown up at $1$ point, as pictured in the toric moment-map diagram below. 

\includegraphics{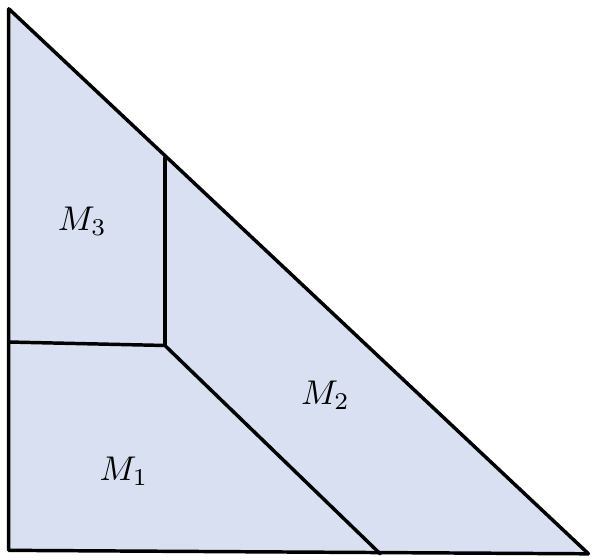}

The picture below is the image of a tropical curve $\gamma$ in the tropical part $\totb{\ex B}$ of  an exploded manifold $\ex B$ with smooth part  fiber pictured above.  The little number at a  vertex in the corner denotes the number of ends of $\gamma$ attached to that vertex, because all  these ends are sent to a point in $\totb {\ex B}$.

\includegraphics{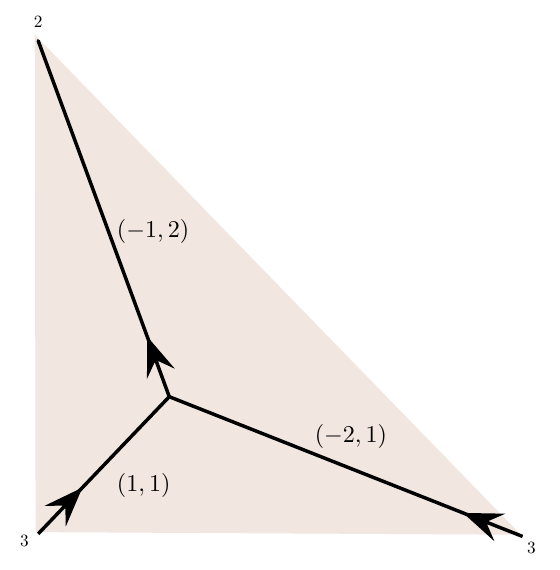}

Let us consider $\eta^{[\gamma_{v}]}$ for the above tropical curve. For $v_{0}$ the central vertex, $\ex B\tc {v_{0}}$ is an exploded manifold $\ex T^{2}$, analogous to $(\mathbb C^{*})^{2}$. The relevant components of $\rend(\ex T^{2})$ for our evaluation map are the quotient of $\ex T^{2}$ by the $\ex T$--actions of weight $(-1,-1)$, $(2,-1)$, and $(-1,2)$ respectively, each isomorphic to $\ex T$. For this vertex, $\eta^{[\gamma_{v_{0}}]}=\hbar\theta_{0}$, where $\theta_{0}$ is the Poincare dual to the product of the three quotient maps $\ex T^{2}\longrightarrow \ex T^{3}$.

For $v_{i}$ the remaining vertices, $\eta^{[\gamma_{v_{i}}]}$ is a Gromov--Witten invariant of curves in the explosion of $M_{i}$ relative to  $M_{i}\cap M_{j}$ and $M_{i}\cap M_{k}$. Curves with tropical part $\gamma_{v_{i}}$ correspond to curves in $M_{i}$ with one special point sent to $M_{1}\cap M_{2}\cap M_{3}$, having contact order $1$ with each divisor, and two or three other special points not contacting any divisor.  Evaluation at these points  is into $\ex T$, and  $\expl M_{i}$ respectively. Our remaining $\eta^{[\gamma_{v_{i}}]}$ are as follows:
\begin{itemize}
\item    $\eta^{[\gamma_{v_{3}}]}=\hbar q^{2E_{31}+E_{32}}$ in the $0$--dimensional cohomology of $(\expl M_{3})^{2}\times \ex T$, where $E_{ij}$ is the symplectic area of $M_{i}\cap M_{j}$.  
\item $\eta^{[\gamma_{v_{1}}]}=\hbar^{2}q^{E_{13}+2E_{12}}\theta_{1}$, where $\theta_{1}\in \rh^{2}((\expl M_{1})^{3}\times \ex T)$ has integral $1$ over the homology class represented by $\ex T$.\footnote{The notation $\rh^{*}$ indicates we are using refined cohomology, from section 9 of \cite{dre}.} 
\item Similar to the case of $v_{1}$, $\eta^{[\gamma_{v_{2}}]}=\hbar^{2}q^{E_{21}+2E_{23}}\theta_{2}$, where $\theta_{2}\in \rh^{2}((\expl M_{2})^{3}\times \ex T)$ has integral $1$ over the homology class represented by $\ex T$.
\end{itemize}

Our gluing formula uses
\[\begin{tikzcd}
\\ (\expl M_{1})^{3}\times(\expl M_{2})^{3}\times (\expl M_{3})^{2}
\\ (\expl M_{1})^{3}\times(\expl M_{2})^{3}\times (\expl M_{3})^{2}\times \ex T^{3}\uar{i^{[\gamma]}}\dar{\Delta}
\\ (\expl M_{1})^{3}\times(\expl M_{2})^{3}\times (\expl M_{3})^{2}\times \ex T^{6}  \end{tikzcd}\]
and reads
\[\eta\tc \gamma=\hbar^{6}q^{2(E_{12}+E_{23}+E_{31})} i_{!}^{[\gamma]}\Delta^{*}(\theta_{1}\wedge \theta_{2}\wedge \theta_{0})\ .\]
As both $\ex T^{3}$ and $\theta_{1}\wedge \theta_{2}\wedge \theta_{0}$ are $6$--dimensional, $\eta\tc\gamma$ is $0$--dimensional, so it suffices to evaluate it at a point. This evaluation amounts to the integral of $\theta_{1}'\wedge\theta_{2}'\wedge \theta_{0}$ over $\ex T^{3}$, where $\theta_{i}'$ is the pullback of the generator of the top-dimensional cohomology on the $i$th  $\ex T$--factor of $\ex T^{3}$. This integral is $3$, so $\eta\tc\gamma=3\hbar^{6}q^{2(E_{12}+E_{23}+E_{31})}$.

 In terms of counting curves, this statement translates to there being $3$ rigid curves in $\ex B$ with tropical part $\gamma$ and the $8$ ends constrained to chosen points in $M_{i}$. Moreover these curves have genus $0$ and symplectic energy $2(E_{12}+E_{23}+E_{31})$, so they are degree $2$. These are not all the curves in $\ex B$ with ends constrained to these points in $M_{i}$  ---  there are $9$ other rigid curves with different tropical parts, as explained in section 8 of \cite{tpgf}. 

\

Let us briefly discuss the Gromov--Witten invariant $\mu$ starring in our second gluing formula, (\ref{gf2}). This invariant uses an enhanced evaluation map
\[EV\co  \Msw_{\cdot}(\ex B)\longrightarrow \mathcal X(\ex B)=\coprod_{g,n}\mathcal X_{g,n}(\ex B)\]
discussed in section \ref{ee map}, where for $2g+-2 +n>0$, 
\[\mathcal X_{g,n}(\ex B)=\Msw_{g,n}(pt)\times_{\cdot/\ex T^{n}}(\End\ex B)^{n}\]
and $\End\ex B$ is the natural stack used to prove gluing formula (\ref{gf1}), and $\Msw_{g,n}(pt)$ is the moduli stack of stable curves mapping to a point.\footnote{$\Msw_{g,n}(pt)$ is constructed in section 4.1 of \cite{evc}. This moduli stack is an exploded orbifold, and the explosion of Deligne-Mumford space relative to its boundary.}
$EV$ is induced from an evaluation map $\ev$ to $(\End\ex B)^{n}$ and a stabilization map $ev^{0}$ to $\Msw_{g,n}$ constructed in section 4.1 of \cite{evc}. In the case that $\totb{\ex B}$ is bounded $\mathcal X_{g,n}=\Msw_{g,n}(pt)\times \ex B^{n}$ and $EV$ is $ev^{0}\times ev$.  In the case that $\ex B$ is the explosion of a manifold with a smooth divisor $D$, the connected component of $\mathcal X_{g,n}$ containing the image of a curve with $n_{2}$ points contacting the divisor, and $n_{1}$ ordinary points is a $\ex T^{n_{2}}$--bundle over $\ex B^{n_{1}}\times \ex D^{n_{2}}\times\mathcal M_{g,n}$.  The precise constructions of $EV$ and  $\mathcal X_{g,n}$ for all $(g,n)$ appear in section \ref{ee map}.

Our Gromov--Witten invariant $\mu$ no longer needs the exponent of $\hbar$ to keep track of Euler class, and is defined  
\[\mu:=EV_{!}(q^{E}) \]
or, in the case of a descendant invariant associated to a complex vectorbundle $W$, 
\[\mu(W):=EV_{!}(q^{E}c(W))\ .\]
The relative invariant $\mu^{[\gamma_{v}]}$ is again the restriction of the corresponding invariant $\mu$ for $\ex B\tc v$ to the connected components 
\[\mathcal X^{[\gamma_{v}]}(\ex B\tc v)\subset \mathcal X(\ex B\tc v)\]
containing the image of curves with tropical part $\gamma_{v}$. Our second gluing formula uses maps
\[\begin{tikzcd}\prod_{v}\mathcal X^{[\gamma_{v}]}(\ex B\tc v)&\lar{\Delta} \ex Y'\rar{I^{[\gamma]}}&\mathcal X(\ex B)\tc{\totb{EV}\gamma}\end{tikzcd}\]
similar to those that appear in our first gluing formula. It reads
\[\mu\tc \gamma=\frac 1{\abs{\Aut\gamma}}I^{[\gamma]}_{!}\Delta^{*}(\prod_{v}\mu^{[\gamma_{v}]})\]
and in the case of descendant invariants using a tautological vectorbundle $W$ on $\Msw(\ex B)$ corresponding to tautological vectorbundles $W_{v}$ on $\Msw_{[\gamma_{v}]}(\ex B\tc v)$, it reads
\[\mu(W)\tc \gamma=\frac 1{\abs{\Aut\gamma}}I^{[\gamma]}_{!}\Delta^{*}(\prod_{v}\mu^{[\gamma_{v}]}(W_{v}))\ .\]

\subsection{Technical assumptions}
Throughout, we shall assume that $\ex B$ is a smooth,  complete, exploded manifold with a taming form $\omega$, and a $\dbar$--log compatible almost-complex structure $J$ tamed by $\omega$.\footnote{For definitions of exploded manifolds and concepts such as basic and complete, see \cite{iec}. Taming forms are discussed in \cite{cem}, section 2, and $\dbar$--log compatible almost complex structures are discussed in section 3, where it is shown that the set of such almost complex structures tamed by $\omega$ is nonempty and contractible.} We shall also assume that the tropical part of $\ex B$ admits an affine immersion into some $\mathbb R^{n}$ so that \cite{cem} establishes the compactness of  the moduli stack of holomorphic curves with bounded energy and genus that map to a connected component of $\coprod_{n}\rend{\ex B}$. This assumption implies that $\ex B$ is basic.

Under these assumptions, \cite{reg,evc} constructs an embedded Kuranishi structure on the moduli stack of curves in $\ex B$. This Kuranishi structure is stronger than many other comparable constructions \cite{FO,MW3,joycebook,pardon} because it comes naturally embedded in a moduli stack  $\Msw(\ex B)$ of smooth, stable, not necessarily holomorphic curves. The virtual fundamental  class of the moduli stack of holomorphic curves is constructed in \cite{vfc}, where we also see how to integrate differential forms over this virtual class, and push forward differential forms over evaluation maps. 

\subsection{Outline of paper}

Our gluing formulae follow from a natural fiber-product diagram 
\begin{equation}\label{cut fp}\begin{tikzcd}[column sep=huge]\Msw_{[\gamma]^{+}}\dar{\cut }\rar& \prod_{e}\bar{\ex B}_{e}\dar{\Delta}
\\ \prod_{v}\Msw_{[\gamma_{v}^{+}]}\rar & \prod_{e}\bar{\ex B}^{2}_{e}\end{tikzcd}\end{equation}
discussed in Theorem \ref{stacks fiber-product}. At the top left of the above diagram is a moduli stack of  curves decorated by a tropical curve $\gamma$, and with an extra choice of point for each  edge of $\gamma$. The downward arrow on the left cuts each such curve at the extra points to obtain cut-curves, and the right-pointing arrows indicate evaluation at those extra points on internal edges, and their cut analogues beneath. This fiber-product diagram, cut-curves and $\gamma$--decorated curves are the subject of section \ref{gluing cut curves}.

Diagram (\ref{cut fp}) does not immediately imply our gluing formulae. Because of the extra choice of points,  the stacks $\Msw_{[\gamma]^{+}}$ and $\Msw_{[\gamma_{v}^{+}]}$ are bundles over the actual stacks $\Msw_{[\gamma]}$ and $\Msw_{[\gamma_{v}]}$ we need for our gluing formulae. To `forget the extra choice of points' we consider a groupoid with objects parametrized by $\Msw_{[\gamma]^{+}}$, and quotient stack $\Msw_{[\gamma]}$. We show, in Proposition \ref{cut*},  that a diagram of groupoids analogous to (\ref{cut fp}) is also fiber-product diagram, however passing to quotient stacks does not quite give a fiber-product diagram. This technical hurdle is overcome using the tool of tropical completion in Lemma \ref{tcfp}, after which our gluing formulae follow quickly from the machinery in \cite{evc,vfc}.

Section \ref{simple ev section} constructs the evaluation map $ev$ relevant for our first gluing formula. To relate to our natural fiber product diagram (\ref{cut fp}), we also construct an evaluation map $\tilde \ev$ on the level of groupoids, inducing a natural evaluation map $\ev$ on the level of quotient stacks. Section \ref{simple gluing formula} contains the proof of our first gluing formula, which is the result of Theorem \ref{sg} and Lemma \ref{eta pf}. In section \ref{ee map}, we enhance our evaluation maps $\ev$ and $\tilde \ev$ to obtain evaluation maps $EV$ and $\tilde{EV}$ that include information about the complex structure of the domain of curves. We then prove our second gluing formula, in Theorem \ref{xg} and Lemma \ref{xpf}. Section \ref{ee map} concludes with an example showing how Kontsevich and Mannin's splitting and genus-reduction axioms follow from our second gluing formula. The  final section of this paper sketches how to extend our gluing formulae to keep track of more topological information, as, for example, is done in \cite{IP} using rim tori.

\section{Gluing cut-curves}
\label{gluing cut curves}

In this section, we define various moduli stacks of (not-necessarily holomorphic) curves, including curves decorated by a tropical curve $\gamma$, and cut-curves. The section concludes with Theorem \ref{stacks fiber-product}, which concerns  a natural fiber-product diagram produced by gluing cut-curves back together. This theorem will be key for proving our tropical gluing formulae. 

 \subsection{ The moduli stack of $\gamma$--decorated curves}

 \

Use the notation $\Msw_{\bullet}\ex B$ for a moduli stack of stable, (not necessarily holomorphic) decorated\footnote{We use decorated moduli stack in the sense of Definition 2.12 of \cite{evc}.} curves in $\ex B$, where the $\bullet$ indicates some unspecified decoration. We shall also have use for possibly unstable curves. Use $\Mpsw_{\bullet}\ex B$ for the moduli stack of (not necessarily stable or holomorphic)  decorated curves in $\ex B$.\footnote{The notation $\Mpsw$  indicates that we are working with exploded manifolds and maps of regularity $\C\infty1$, which is as good as smooth for all practical purposes.}  We shall use $\Mu_{\bullet}$ in statements that hold for both $\Msw_{\bullet}$ and $\Mpsw_{\bullet}$.

Let $\gamma$ be a tropical curve in $\totb{\ex B}$. Let $\Mu_{\gamma}$  indicate the moduli stack of curves within $\Mu$   with a chosen isomorphism of their tropical part to $\gamma$.  Below we shall define a decorated moduli stack $\Mu_{[\gamma]}$ that is a kind of closure of $\Mu_{\gamma}$.  Loosely, $\Mu_{[\gamma]}$ is the closure of $\Mu_{\gamma}$ in the moduli stack of  curves $f$ in $\ex B$ with tropical part suitably labeled by $\gamma$.

 \begin{defn}[$\gamma$--decoration]\label{decorated}
 A $\gamma$--decorated tropical curve is  a tropical curve in $\totb{\ex B}$ with a continuous affine map of its domain to the domain of $\gamma$ so that this map  is a homeomorphism restricted to the inverse image of the interior of all edges of $\gamma$ and  is an integral-affine isomorphism restricted to all exterior edges.  
 
 Define a $\gamma$--decorated curve to be a curve with a $\gamma$--decorated tropical part.  Consider the stack  of $\gamma$--decorated curves, with objects   $\C\infty1$ families of curves $\hat f\in \Mu_{\cdot}(\ex B)$ with an affine map of $\totb{\ex C(\hat f)}$ to the domain of $\gamma$ that makes each individual curve a $\gamma$--decorated curve.   $\Mu_{\gamma}$ is a substack of this stack of $\gamma$--decorated curves. Define  $\Mu_{[\gamma]}$ to be the closure of $\Mu_{\gamma}$ in this 
stack of $\gamma$--decorated curves.  

\end{defn}

%

To prove our gluing formula, we consider the following bundle over $\Mu_{[\gamma]}$. 

\begin{defn}\label{stack+}Define $\Mu_{[\gamma]^{+}}$ to be the stack of curves in $\Mu_{[\gamma]}$  with an additional choice of point in each edge labeled by an  edge of $\gamma$. So, a family in $\Mu_{[\gamma]^{+}}$ is a family $\hat f$ in $\Mu_{[\gamma]}$, along with, for each  edge $e$ of $\gamma$, a $\C\infty1$ section 
\[s_{e}\co\ex F(\hat f)\longrightarrow \ex C( \hat f)\]
 with image contained in the edge labeled by $e$.

Composing $s_{e}$ with $\hat f$ defines a $\C\infty1$ evaluation map
\[ev_{e}\co\Mu_{[\gamma]^{+}}\longrightarrow  \ex B\]
 \end{defn}
Let $\bar{\ex B}_{e}$ indicate the closure of the strata of $\ex B$ with tropical part containing the edge $e$. The evaluation map $ev_{e}$ always lands in $\bar{\ex B}_{e}$.

\[ev_{e}\co\Mu_{[\gamma]^{+}}\longrightarrow  \bar{\ex B}_{e}\]

  This stack $\Mu_{[\gamma]^{+}}$ is a  bundle over $\Mu_{[\gamma]}$. In particular, given any family $\hat f$ in $\Mu_{[\gamma]}$, the pullback of $\hat f$ to $\Mu_{[\gamma]^{+}}$ is a family $\hat f^{+}$ constructed as follows: Let $\ex C_{e}(\hat f)$ indicate the strata of $\ex C(\hat f)$ decorated by $e$. The domain $\ex F(\hat f^{+})$ is the fiber-product of all these $\ex C_{e}(\hat f)$ over $\ex F(\hat f)$, the family $\hat f^{+}$ is the pullback of $\hat f$ under the map $\ex F(\hat f^{+})\longrightarrow \ex F(\hat f)$, and $s_{e}$ is given by the canonical map $\ex F(\hat f^{+})\longrightarrow \ex C_{e}(\hat f)$.
   
  The family $\hat f^{+}$ is the pullback of $\hat f$ to $\Mu_{[\gamma]^{+}}$ in the sense that it satisfies the following universal condition. Let $\pi\co\Mu_{[\gamma]^{+}}\longrightarrow \Mu_{[\gamma]}$ indicate the map that forgets the extra sections $s_{e}$.  Given any family $\hat g$ in $\Mu_{[\gamma]^{+}}$ and a map $\pi(\hat g)\longrightarrow \hat f$, there exists a unique lift of this map making the diagram below commute.\footnote{ Unlike the other arrows in the above diagram,  the arrows $\pi$ are not maps of families of curves.}
  \begin{equation}\label{f+ universal property}\begin{tikzcd}[row sep= small]\hat g \rar[dotted]{\exists !}\dar{\pi}& \hat f^{+}\dar{\pi}
\\ \pi(\hat g)\ar{dr}\rar[dotted]&\pi(\hat f^{+})\dar 
\\ & \hat f\end{tikzcd}\end{equation}

Both $\ex C_{e}(\hat f)$ and $\bar{\ex B}_{e}$ have a (not everywhere defined) action of $\ex T$ so that $ev_{e}$ is equivariant.  In section \ref{simple ev section}, we   encode this action as a groupoid structure on $\ex C_{e}(\hat f)$ and $\bar{\ex B}_{e}$, and extend $ev_{e}$ to a  map of groupoids. Taking  quotients will then define an evaluation map from $\Mu_{[\gamma]}$ instead of $\Mu_{[\gamma]^{+}}$.

  \subsection{Cut-curves}

  The goal of this section is to prove that $\Mu_{[\gamma]^{+}}$ is a fiber-product of some moduli stacks of cut-curves. For this, we shall be assuming that the domain of  $\gamma$ is connected, and  not $\mathbb R$  ---  in this exceptional case, $\Mu_{[\gamma]}$ is easily studied directly. 
    We can cut any curve in $\Mu_{[\gamma]^{+}}$ at the extra points on its edges to obtain cut-curves.     Cut-curves were defined in \cite{elc} omitting the definition of a family of cut-curves, given below.

\begin{defn}\label{cut family} A family of cut-curves over an exploded manifold $\ex F$ is 
\begin{enumerate}
\item an abstract exploded space\footnote{ See \cite{iec}, Definition 3.1.  An abstract exploded space is a topological space with a sheaf  of $\mathbb C^{*}\e{\mathbb R}$--valued functions. One way of defining an abstract exploded space is as a subset of an exploded manifold $\ex B$  given the subspace topology and the pullback of the sheaf of functions $\mathcal E^{\times}(\ex B)$.} $\ex C$ with a map of abstract exploded spaces $\ex C\longrightarrow \ex F$, 
\item some number of sections $s_{e}\co\ex F\longrightarrow \ex C$ called cuts so that $\ex C$ minus the image of these sections is an exploded manifold,
\item a fiberwise almost complex structure $j$ on $\ex C$ minus all cuts,
\end{enumerate}
so that the above data is the result of trimming  some family of curves $ (\ex C',j)$ over $\ex F$ at sections $s_{e}\co\ex F\longrightarrow \ex C'$ with image in distinct ends of curves in $\ex C'$; here `trimming' means that the above data is the restriction of $(\ex C',j,\{s_{e}\})$ to the subset of $\ex C'$ consisting of the image of the sections $s_{e}$ and all points in $\ex C'$ that  have tropical part not as far out on any edge  as the tropical part of the image of $s_{e}$.

A family of (holomorphic) cut-curves in $\ex B$ is a (fiberwise holomorphic) map $\ex C\longrightarrow \ex B$.

\end{defn}

Although the domain of a cut-curve  is the result of trimming some honest curve, it is not true that every cut-curve in $\ex B$ can be obtained by trimming a curve in $\ex B$; this is the case if and only if all the cut edges of the corresponding tropical curve  extend  semi-infinitely  in $\totb{\ex B}$.

Each cut-curve $f$ in $\ex B$ has a tropical part $\totb{f}$ with special $1$--valent vertices at the cuts. Call such a tropical curve a cut tropical curve. We may define a moduli space $\Mu_{[\gamma]}$ of $\gamma$--decorated cut-curves as in Definition \ref{decorated}.

Given a tropical curve $\gamma$ with a choice of point on each edge, we can obtain a (possibly disconnected) cut tropical curve by cutting $\gamma$ at these chosen points, and discarding\footnote{In our neglected exceptional case that the domain of $\gamma$ is $\mathbb R$, this cutting process would discard everything.} cut edges not attached to a vertex $v$. This cut-curve will have one connected component $\gamma_{v}^{+}$ for each vertex $v$ of $\gamma$, so we can write this cut-curve as $\coprod_{v}\gamma_{v}^{+}$. 

\

\begin{lemma}\label{cut lemma}There exists a map of stacks
\[\cut\co \Mu_{[\gamma]^{+}}\longrightarrow \prod_{v}\Mu_{[\gamma_{v}^{+}]}\ .\]
In particular, given a family of curves $\hat f$ in $\Mu_{[\gamma]^{+}}$, there exists a family of curves $\cut \hat f$ in $\prod_{v}\Mu_{[\gamma^{+}_{v}]}$ and a morphism $\ex C(\cut\hat f)\longrightarrow \ex C(\hat f)$ so that the following diagram commutes
\[\begin{tikzcd}\ex C(\cut\hat f)\ar[bend left]{rr}{\cut\hat f}\rar\dar&\ex C(\hat f)\dar\rar[swap]{\hat f}& \ex B
\\ \ex F(\cut\hat f)\rar{\cong}\uar[bend left]{s_{e_{i}}}&\ex F(\hat f)\uar[bend left]{s_{e}}\end{tikzcd}\]
where  each internal edge $e$ of $\gamma$ corresponds to two cut edges, $e_{1}$ and $e_{2}$, of $\coprod_{v}\gamma_{v}$, (and an end  of $\gamma$ corresponds to a unique cut edge of $\coprod_{v}\gamma_{v}$). Moreover, any morphism $x:\ex C(\cut\hat f)\longrightarrow \ex A$ is the pullback of a morphism $x':\ex C(\hat f)\longrightarrow \ex A$ if and only if for all internal edges $e$ of $\gamma$,  $x\circ e_{1}=x\circ e_{2}$ and the derivative of $\totb x$ along the edges $e_{i}$ is opposite,  and for all ends of $\gamma$, $\totb{x}$ is infinitely extendible along the corresponding edge of $\coprod_{v}\gamma_{v}$. 
\end{lemma}

\pf

The idea is to cut our family  $\hat f$ at the sections $s_{e}\co\ex F(\hat f)\longrightarrow \ex C(\hat f)$ to obtain a family of  cut-curves (still parametrized by  $\ex F(\hat f)$) with one connected component $\hat f_{v}$ in $\Mu_{[\gamma_{v}^{+}]}$ for each vertex $v$ of $\gamma$.

If there are no loops  attached to $v$, we can define $\ex C(\hat f_{v})$  as a subset of $\ex C(\hat f)$ with the induced exploded structure. Namely,  $\ex C(\hat f_{v})$ is the union of the  image of the sections $s_{e}\co\ex F(\hat f)\longrightarrow \ex C(\hat f)$ for all edges $e$ adjacent to $v$, and  all points in $\ex C$ with tropical part fiberwise closer to $v$ than these sections. It is easily checked that $\ex C(\hat f_{v})\longrightarrow \ex F(\hat f)$ along with these sections $s_{e}$, and the restriction of $\hat f$ to $\ex C(\hat f_{v})$ is a family of cut-curves. The above procedure fails in the case that there is an edge $e$ with two ends attached to $v$,  as we  need to break apart the two halves of this edge. 
 
 In the general case, define the set $\ex C(\hat f_{v})$ as the union of a copy of  $\ex F$ for each oriented edge leaving $v$ with the subset of $\ex C(\hat f)$ with tropical part sent closer to $v$ than $s_{e}$. $\ex C(\hat f_{v})$ comes with a natural map to $\ex C(\hat f)$ which is the identity inclusion on the main part, and $s_{e}$ on each copy of $\ex F$ corresponding to an edge $e$.
 
  Before continuing with the description of $\ex C(\hat f_{v})$, consider the case of an individual curve $f$ in $\hat f$. The $\gamma$--decoration gives a map $\ex C(f)\longrightarrow \gamma$, and the image of $s_{e}$ under this map is sent to a point on the edge $e$. Cut $\gamma$ at these  points, and consider the connected component $\gamma_{v}^{+}$ containing $v$. As a set, $\ex C(f_{v})$ is  the inverse image of the interior of $\gamma_{v}^{+}$ with a point for each oriented edge leaving $v$; the identity inclusion together with the maps $s_{e}$ defines a natural map $\ex C( f_{v})\longrightarrow \ex C(f)$ that is injective everywhere apart from the extra points corresponding to an edge $e$ with both ends attached to $v$. With this understood, we could equivalently define $\ex C(\hat f_{v})$ as the union of $\ex C(f_{v})$ for all $f$ in $v$. Note that $\ex C(f_{v})$ comes with a surjective map to $\gamma_{v}^{+}$ so that the following diagram commutes.
 \[\begin{tikzcd}\ex C(f_{v})\dar\rar &\ex C(f)\dar
 \\ \gamma^{+}_{v}\rar &\gamma\end{tikzcd}\]
The point corresponding to an oriented edge leaving $e$ is sent to the cut  endpoint of that edge in $\gamma_{v}^{+}$, and elsewhere, the map $\ex C(f_{v})\longrightarrow \gamma^{+}$ is uniquely determined by the above diagram. The topology on $\ex C(f_{v})$ is the topology generated by the inverse image of open sets from $\ex C(f)$, and  the inverse image of the closure of any edge of $\gamma_{v}^{+}$.
 
 We can describe the topology on $\ex C(\hat f_{v})$ similarly. For other curves $f$ in $\hat f$, the corresponding $\gamma^{+}_{v}$ may have different edge lengths, however there is a canonical piecewise-linear isomorphism between all these $\gamma^{+}_{v}$, so we obtain a $\gamma_{v}^{+}$--decoration map $\ex C(\hat f_{v})\longrightarrow \gamma_{v}^{+}$. The topology on $\ex C(\hat f_{v})$ is the topology generated by the inverse image of open subsets of $\ex C(\hat f)$ and the inverse image of the closure of any edge in $\gamma_{v}^{+}$.
 With this topology, we can define the sheaf of exploded functions as follows: for any open subset $U$ small enough that its closure only intersects (the inverse image of) one edge of $\gamma_{v}^{+}$, the exploded functions are the $\mathbb C^{*}\e{\mathbb R}$--valued functions pulled back from exploded functions from $\ex C(\hat f)$ defined some open neighborhood of the image of $U$.
 
 Note that as exploded functions on $\ex C(\hat f_{v})$ are locally the pullback of exploded functions on $\ex C$,  any morphism $x:\ex C(\cut\hat f)\longrightarrow \ex A$ is the pullback of a morphism $x':\ex C(\hat f)\longrightarrow \ex A$ if and only if for all internal edges $e$ of $\gamma$,  $x\circ e_{1}=x\circ e_{2}$ and the derivative of $\totb x$ along the edges $e_{i}$ is opposite,  and for all ends of $\gamma$, $\totb{x}$ is infinitely extendible along the corresponding edge of $\coprod_{v}\gamma_{v}$.

 With this definition, it is clear that $\ex C(\hat f_{v})\longrightarrow \ex C(\hat f)$ is a morphism of abstract exploded manifolds, and that this construction is functorial: given any morphism $\hat g\longrightarrow \hat f$ in $\Mu_{[\gamma]}$, there is a unique morphism $\ex C(\hat g_{v})\longrightarrow \ex C(\hat f_{v})$ so that the following diagram commutes.
 \[\begin{tikzcd}\ex C(\hat g_{v})\dar\rar&\ex C(\hat f_{v})\dar
 \\ \ex C(\hat g)\rar& \ex C(\hat f)\end{tikzcd}\]
For each oriented edge $e$ of $\gamma_{v}^{+}$, there is a unique map (of sets)\[s_{e}\co\ex F(\hat f_{v}):=\ex F(\hat f)\longrightarrow \ex C(\hat f_{v})\] that is a section of $\ex C(\hat f_{v})\longrightarrow \ex F(\hat f_{v})$ with image in the inverse image of the cut end of $e$ in $\gamma_{v}^{+}$. The fact that the diagram below commutes implies that $s_{e}$ is a morphism of abstract exploded manifolds; these morphisms $s_{e}$ define the cuts of our family of cut-curves. 
\[\begin{tikzcd}\ex C(\hat f_{v})\rar&\ex C(\hat f)
\\ \ex F(\hat f_{v})\rar\uar{s_{e}}&\ex F(\hat f)\uar{s_{e}}\end{tikzcd}\]
  The rest of the data for defining $\hat f_{v}$ as a family of cut-curves is the  fiberwise almost complex structure pulled back from $\ex C(\hat f)$, and the map $\hat f_{v}\co\ex C(\hat f_{v})\longrightarrow \ex B$ that is the pullback of $\hat f\co\ex C(\hat f)\longrightarrow \ex B$. With this definition it is easy to verify that given any morphism $\hat g\longrightarrow \hat f$ in $\Mu_{[\gamma]}$, the unique maps $\ex C(\hat g_{v})\longrightarrow \ex C(\hat f_{v})$ above indeed define a unique morphism $\cut\hat g\longrightarrow \cut\hat f$ compatible with all this structure, so our construction is functorial. With the functoriality of the construction understood, it is also easy to verify locally that the resulting $\hat f_{v}$ is indeed a family of cut-curves in $\Mu_{[\gamma^{+}]}$ satisfying Definition \ref{cut family}.
  
  We have therefore defined a map of stacks, 
  \[\cut \co\Mu_{[\gamma]^{+}}(\ex B)\longrightarrow \prod_{v}\Mu_{[\gamma_{v}^{+}]}(\ex B)\]
  where $\cut (\hat f):=\coprod_{v}\hat f_{v}$, the (possibly disconnected) family of cut-curves  parametrized by $\ex F(\hat f)$ with domain $\coprod_{v}\ex C(\hat f_{v})$. 
  
  \stop

%
%
%

\

For the following theorem, use the notation $\ie \gamma$ for the set of internal edges of $\gamma$. For each $e\in \ie\gamma$, there are two corresponding  cut edges $e_{1}, e_{2}$ of $\coprod_{v}\gamma_{v}^{+}$ where we can define evaluation maps $ev_{e_{i}}$. 

\begin{thm}\label{stacks fiber-product}The following is a fiber-product diagram.
\[\begin{tikzcd}[column sep=huge]\Mu_{[\gamma]^{+}}\dar{\cut }\ar{rr}{\prod_{e\in\ie\gamma} ev_{e}}&& \prod_{e}\bar{\ex B}_{e}\dar{\Delta}
\\ \prod_{v}\Mu_{[\gamma_{v}^{+}]}\ar{rr}{\prod_{e\in\ie\gamma}(ev_{e_{1}},ev_{e_{2}})} && \prod_{e}\bar{\ex B}^{2}_{e}\end{tikzcd}\]

More precisely,  given any  family of cut-curves $\hat f$ in $\prod_{v}\Mu_{[\gamma_{v}^{+}]}$ so that for all internal edges $e$ of $\gamma$, $ev_{e_{1}}=ev_{e_{2}}$ on $\hat f$, there exists a family $\cut^{*}\hat f$ in $\Mu_{[\gamma]^{+}}$ with an isomorphism $\cut (\cut^{*}\hat f)\longrightarrow \hat f$. This family $\cut^{*}\hat f$ satisfies the universal property that  given any other family of curves $\hat g$ in $\Mu_{[\gamma]^{+}}$ with a map $\cut \hat g\longrightarrow \hat f$, there exists a unique map $\hat g\longrightarrow \cut^{*}\hat f  $ so that the following diagram commutes:

 \[\begin{tikzcd}\hat g \rar[dotted]{\exists !}\dar{\cut }& \cut^{*}\hat f  \dar{\cut }
\\ \cut \hat g\ar{dr}\rar[dotted]&\cut (\cut^{*}\hat f  )\dar{\simeq} 
\\ & \hat f\end{tikzcd}\]

\end{thm}

\pf

%

The family $\cut^{*}\hat f  $ is constructed by gluing together $\hat f$ at the matching cuts $s_{e_{1}}$ and $s_{e_{2}}$. One way to describe $\ex C(\cut^{*}\hat f  )$ is as follows:
\begin{itemize}
\item The set of points in $\ex C(\cut^{*}\hat f  )$ is the union of $\ex C(\hat f)\setminus(\bigcup_{e_{i}} s_{e_{i}}(\ex F(\hat f)))$ with $\mathbb C^{*}\times \ex F(\hat f)$ for each internal edge $e$, and $\mathbb C^{*}\e{[0,\infty)}\times \ex F(\hat f)$ for each end of $\gamma$. Write these extra points $(c,p)$ as $c *s_{e}(p)$. 

There is a canonical map $\ex C(\hat f)\longrightarrow \ex C(\cut^{*}\hat f  )$: this map is  the identity on $\ex C(\cut^{*}\hat f  )$ minus all cuts, and $s_{e_{i}}(p)\mapsto 1*s_{e}(p)$ on cuts.
\item $\ex C(\cut^{*}\hat f  )$ has the following topology: The open subsets $U\subset \ex C(\cut^{*}\hat f  )$ have  open inverse image in $\ex C(\hat f)$ and  satisfy the additional conditions  that $c*s_{e}(p)\in U$ if and only if $1*s_{e}(p)\in U$, and $1*s_{e_{1}}(p)\in U$ if and only if $1*s_{e_{2}}(p)\in U$.

\item The defining sheaf of exploded functions on $\ex C(\cut^{*}\hat f  )$ is  as follows. The exploded functions $x$ on $U\subset \ex C(\cut^{*}\hat f)$ are those that pull back to exploded functions on $\ex C(\hat f)$, and satisfy the following addition conditions
\begin{itemize}\item $x\circ s_{e_{1}}=x\circ s_{e_{2}}$ 
\item The derivative  of $\totb{x}$ on the edges $e_{1}$ and $e_{2}$ is opposite ---  in other words, if $\tilde z_{i}$ indicates standard coordinates on these edges, then $\tilde z_{1}\frac \partial {\partial z_{1}}x=\alpha x$ and $\tilde z_{2}\frac \partial {\partial z_{2}}x=-\alpha x$.
\item  If $\alpha$ indicates the derivative of $\totb x$ on the edge $e$, then \[x(c*s_{e}(p)):=c^{\alpha}x(s_{e_{i}}(p))\ .\]
\end{itemize}
The above implies that maps $x$ from $\ex C(\cut^{*}\hat f  )$ to any exploded manifold canonically correspond to maps $x$ from $\ex C(\hat f)$ so that $x\circ s_{e_{1}}=x\circ s_{e_{2}}$, the derivative of $\totb{x}$ on the edges $e_{1}$ and $e_{2}$ is opposite, and the remaining external edges  of tropical curves in $\totb x$ can be extended to be semi-infinite. In particular, the map $\ex C(\hat f)\longrightarrow \ex F(\hat f)$ induces a map $\ex C(\cut^{*}\hat f  )\longrightarrow \ex F(\hat f)=\ex F(\cut^{*}\hat f  )$. It also follows that the section $s_{e}\co\ex F(\cut^{*}\hat f  )\longrightarrow \ex C(\cut^{*}\hat f  )$ pulls back exploded functions to exploded functions, so is a valid map of exploded manifolds.
\end{itemize}

To see that the above defines an exploded manifold structure on $\ex C(\cut^{*}\hat f  )$, it suffices to check locally around the section $s_{e}$. We do this for $e$ an internal edge, the argument for an external edge is similar but easier. Around any point in $\ex F(\hat f)$ there exists a coordinate chart $U$  so that a neighborhood of $s_{e_{i}}(U)$ is isomorphic to the result of trimming $U\times \et 11$ by sections $s_{e_{i}}\co U\longrightarrow \et 11$. Then, there is an open  neighborhood   $V$ of  $s_{e}(U)\subset \ex C(\cut^{*}\hat f  )$ with topology and exploded functions the same as the following fiber-product:
\[\begin{tikzcd} V\rar\dar& U\times \et 11\times \et 11\dar{(u,\tilde z_{1},\tilde z_{2})\to (u,\tilde z_{1}\tilde z_{2})}
\\ U\rar{s_{e_{1}}s_{e_{2}}} & U\times \et 11\end{tikzcd}\]
Therefore,  $\ex C(\cut^{*}\hat f  )$ as defined above really is an exploded manifold. The fiberwise almost complex structure defining  $\ex C(\cut^{*}\hat f  )\longrightarrow \ex F(\cut^{*}\hat f  )$ as a family of  curves is induced from  $\ex C(\hat f)\longrightarrow \ex F(\hat f)$ so that the inclusion $\ex C(\hat f)\longrightarrow \ex C(\cut^{*}\hat f)$ is fiberwise holomorphic. 

We have already assumed that  $\hat f\co  \ex C(\hat f)\longrightarrow \ex B$ satisfies the condition that $\hat f\circ s_{e_{1}}=\hat f\circ s_{e_{2}}$, so Lemma \ref{cut lemma} implies that to verify that $\hat f$ induces a map $\ex C(\cut^{*}\hat f  )\longrightarrow \ex B$, it suffices to check that the derivatives of $\totb{\hat f}$ along the  edges $e_{1}$ and $e_{2}$ are opposite, and that the remaining external edges of the tropical curves in $\totb f$ extend  semi-infinitely. These conditions hold if $\hat f$ is in $\prod_{v}\Mu_{\gamma_{v}^{+}}$, and  therefore hold for $\hat f$ in $\prod_{v}\Mu_{[\gamma_{v}^{+}]}$, because $\prod_{v}\Mu_{\gamma_{v}^{+}}$ is dense in $\prod_{v}\Mu_{[\gamma_{v}^{+}]}$. Our map $\hat f$ therefore induces our map $\cut^{*}\hat f$. 
\[\begin{tikzcd}\ex C(\hat f)\ar[bend right]{rr}{\hat f}\rar&\ex C(\cut^{*}\hat f)\rar{\cut^{*}\hat f}& \ex B\end{tikzcd}\]
Note that $\cut^{*}\hat f$ consists of stable curves if and only if $\hat f$ consists of stable curves.
 There is a unique isomorphism $\cut (\cut^{*}\hat f  )\longrightarrow \hat f$ compatible with our two canonical maps $\ex C(\hat f)\longrightarrow \ex C(\cut^{*}\hat f  )$
and $\ex C(\cut (\cut^{*}\hat f  ))\longrightarrow \ex C(\cut^{*}\hat f  )$. 

Now we check the claimed universal property of $\cut^{*}\hat f  $. Given a family $\hat g$ with a map $\cut \hat g\longrightarrow \hat f $ in $\prod_{v}\Mu_{[\gamma_{v}^{+}]}$, the composition  $\ex C(\cut \hat g)\longrightarrow \ex C(\hat f)\longrightarrow \ex C(\cut^{*}\hat f  ) $ satisfies the conditions from Lemma \ref{cut lemma} to be the pullback of a map $\ex C(\hat g)\longrightarrow \ex C(\cut^{*}\hat f)$. This map defines a  unique map $\hat g\longrightarrow \cut^{*}\hat f$ so that the required diagram commutes.

\[\begin{tikzcd}[row sep= small]\hat g \rar[dotted]{\exists !}\dar{\cut }& \cut^{*}\hat f  \dar{\cut }
\\ \cut \hat g\ar{dr}\rar[dotted]&\cut (\cut^{*}\hat f  )\dar{\simeq} 
\\ & \hat f\end{tikzcd}\]   
 
\stop 

The fiber-product diagram from Theorem \ref{stacks fiber-product} is equivariant with respect to some (partially defined) actions of $\ex T$ corresponding to the extra choices of points on curves in $\Mu_{[\gamma]^{+}}$ and the length of cut edges in $\Mu_{[\gamma_{v}^{+}]}$. We shall formalize this in the language of groupoids  in the next section. For our gluing formula, and  to define  evaluation maps from $\Mu_{[\gamma]}$ instead of $\Mu_{[\gamma]^{+}}$, we must take the quotient by these $\ex T$ actions. 

\section{The simple evaluation map}\label{simple ev section}

%

Each end of a curve $\ex C\neq \ex T$ is isomorphic to $\et 1{(0,\infty)}$, however there is a $\mathbb C^{*}$-fold choice of isomorphism. The moduli stack of maps $\et 1{(0,\infty)}\longrightarrow \ex B$, identified up to isomorphism, is not a nice stack, so we shall replace it with a better behaved stack for defining our simple evaluation map. In particular, we  throw away  `tropical' information  by  identifying two maps $\et1{(0,\infty)}\longrightarrow \ex B$ if they eventually coincide.  By doing so, we  obtain a better behaved stack $\End\ex B$, which is a quotient of an exploded  manifold $\rend\ex B$ by some trivial group actions, and is an orbifold on most components. 

  In what follows, we shall describe $\End\ex B$ as the quotient stack of a Lie groupoid $\gend\ex B$ consisting of a collection of connected components of a Lie groupoid $\mathcal G(\ex B)$ with objects parametrized by the integral-vectors ${}^{\mathbb Z}T\ex B$.

We shall need various partially defined $\ex T$--actions, all in the form a `flow' induced by an integral vector. Recall, from Definition 6.8 of \cite{iec}, that integral tangent vectors ${}^{\mathbb Z}T\ex B\subset T\ex B$ are  vectors $v$ so that $v\tilde z$ is an integer times $\tilde z$ for every (locally defined) exploded function $\tilde z$. For example, in standard coordinates on $\et 1{[0,\infty)}$, the integral-vectors consist of the zero-vector, and integer multiples of the real part of $\tilde z\frac\partial{\partial\tilde z}$ on the strata over $(0,\infty)$, where $\totb{\tilde z}=0$. As a subset of the exploded manifold $T\ex B$, the set of integral tangent vectors inherits the structure of an exploded manifold. For example ${}^{\mathbb Z}T\et 1{[0,\infty)}$ is isomorphic to the disjoint union of  $\et 1{[0,\infty)}$ with a copy of $\et 1{(0,\infty)}$ for every nonzero integer. There is a natural (not everywhere defined) action of $\ex T$ on ${}^{\mathbb Z}T\ex B$ so that for any  exploded function $\tilde z$, 
\[(c*v)\tilde z=c^{\tilde z^{-1}v\tilde z}v\tilde z\ .\]
There may be no vector $c*v$ satisfying the above conditions, in which case $c*v$ is not defined. For example, if $v$ is the real part of $\tilde z\frac\partial{\partial \tilde z}$ over the point where $\tilde z=1\e{1}$, $c*v$ is  the same vector at the point $\tilde z=c\times\e{1}$ so long as $\totb c>-1$. For $\totb c\leq -1$,  $c^{*}v$ is not defined. For $c\in \mathbb C^{*}\subset \ex T$, this action is always defined, and $e^{t}*v$ can be thought of as  the result of flowing $v$ for time $t$.

The above action defines a groupoid $\mathcal G(\ex B)$ with objects $\mathcal G_{0}(\ex B)$ parametrized by  ${}^{\mathbb Z}T\ex B$, and morphisms $\mathcal G_{1}(\ex B)$ parametrized by the set of $(c,v)$ so that $c*v$ is defined.  We shall see that $\mathcal G_{1}(\ex B)$ inherits an exploded manifold structure as a subset of $\ex T\times \mathcal G_{0}(\ex B)$.  This groupoid is a nice Lie groupoid in the category of exploded manifolds, with  all structure maps morphisms in the category of exploded manifolds. Moreover, $\mathcal G$ defines a functor from the category of exploded manifolds to the category of such Lie groupoids. 

To understand this functor $\mathcal G$, consider $\mathcal G(\et mP)$. For each integral-vectorfield $v$ on $P$, there is a corresponding integral-vectorfield on the subset of $\et mP$ with tropical part $P_{v}$ the union of all strata of $P$ tangent to $v$.  Then 
\[\mathcal G_{0}(\et mP)=\coprod_{v\in\mathbb Z^{m}}\et m{P_{v}}\ .\]
Morphisms in $\mathcal G(\et mP)$ always stay within these connected components, and $\mathcal G_{1}(\et mP)$ is also a union of connected components indexed by integral-vectors $v$. Define $\hat P_{v}\subset P_{v}\times \mathbb R$ to be
\[\hat P_{v}:=\{(p,x)\in P_{v}\times \mathbb R\text{ so that }p+xv\in P_{v}\}\ .\]
The tropical part of $\mathcal G_{1}(\et mP)$ is the disjoint union of these $\hat P_{v}$.
\[\mathcal G_{1}(\et mP)=\coprod_{v\in\mathbb Z^{m}}\et {m+1}{\hat P_{v}}\]
There is a canonical inclusion of $P_{v}$ into $\hat P_{v}$ as $P_{v}\times 0$, corresponding to the `identity' section $\id\co  \et m{P_{v}}\longrightarrow \et {m+1}{\hat P_{v}}$ sending $v\in \et m{P_{v}}$ to $(v,1\e 1)$.  $\hat P_{v}$ comes with two surjective integral-affine maps to $P_{v}$: the projection $(p,x)\mapsto p$, and the map $(p,x)\mapsto p+xv$. These two maps are the tropical part of the two structure maps $s,t\co  \mathcal G_{1}\longrightarrow \mathcal G_{0}$. 
\[s,t\co  \et {m+1}{\hat P_{v}}\longrightarrow \et m{P_{v}} \]
\[s(v,c):=v\]
\[t(v,c):=c*v\]
Put the analogous exploded manifold structure on  $\mathcal G(\mathbb R^{n}\times \et mP)=\mathbb R^{n}\times \mathcal G(\et mP)$. To put an exploded manifold structure on $\mathcal G_{1}(\ex B)$, note that there is a natural projection $\mathcal G_{0}(\ex B)\longrightarrow \ex B$, and that our action is always trivial on $\totl{{}^{\mathbb Z}T\ex B}$, so $s$ and   $t$ followed by the smooth part of this projection define the same map $\mathcal G_{1}(\ex B)\longrightarrow \totl{\ex B}$. Put the pulled-back topology on $\mathcal G_{1}(\ex B)$, and then give $\mathcal G(\ex B)$ over a coordinate chart  the exploded structure discussed above. This exploded structure on $\mathcal G(\ex B)$ is well-defined, and  $\mathcal G$ defines a functor to the category of (exploded) Lie groupoids.  

\begin{remark}\label{Gfp} The functor $\mathcal G$  sends fiber-products to fiber-products. 
\[\mathcal G(\ex A\times _{\ex C}\ex B)=\mathcal G (\ex A)\times_{\mathcal G(\ex C)}\mathcal G(\ex B)\]
The special case of $\ex A\times \ex B$ must be considered as the fiber-product over a point; $\mathcal G(pt)$ is the groupoid with a single object, and morphisms parametrized by $\ex T$.
\end{remark}

If we consider $\mathcal G(\ex B)$ as a stack $\mathcal{BG}(\ex B)$ (in other words, replace $\mathcal G(\ex B)$ by the quotient stack classifying principle $\mathcal G(\ex B)$--bundles, as in Definition 3.17 of \cite{orbifoldstack} translated to the category of exploded manifolds), then the components of $\mathcal {BG}$ corresponding to nonzero integral-vectors are  orbifolds,\footnote{By an exploded orbifold, we mean a Deligne-Mumford stack in the category of exploded manifolds} and the components corresponding to primitive integral-vectors are  exploded manifolds. The remaining zero-vector component is  the quotient of $\ex B$ by the trivial $\ex T$--action.

More explicitly, $\mathcal {BG}$ consists of the stack of bundles $\ex L\longrightarrow \ex F$ along with maps of groupoids,
\[\begin{tikzcd}\ex L\times_{\ex F}\ex L\rar{\phi_{1}}\dar{\pi_{1}\text{ or }\pi_{2}}&\mathcal G_{1}\dar{s\text{ or }t}
\\ \ex L\rar{\phi_{0}}&\mathcal G_{0}\end{tikzcd}\]
so that the above diagrams are fiber-product diagrams.  The above data defines an object in $\mathcal {BG}$, and a morphism is a commutative diagram
\[\begin{tikzcd}\ex L_{1}\rar\dar&\ex L_{2}\dar
\\ \ex F_{1}\rar&\ex F_{2}\end{tikzcd}\]
compatible with the maps into $\mathcal G$.

For example, consider the component of $\mathcal {BG}(\et mP)$ corresponding to a nonzero integral-vector $v$. The  $v$--component of $\mathcal G_{0}(\et mP)$ is $\et m{P_{v}}$. Each fiber of $\ex L\longrightarrow \ex F$ must be isomorphic to some $\et 1{(a,b)}$, and the map $\phi_{0}$ to $\et m{P_{v}}$ restricted to each fiber must have tropical part with derivative $v$. Moreover, $\phi_{0}$ restricted to each fiber must be complete. (The map $\phi_{1}$ is uniquely determined by  $\phi_{0}$.) Let $P/v$ be the image of $P_{v}$ under the projection with kernel spanned by $v$.    There is a unique map $\ex F\longrightarrow \et {m-1}{P/v}$ so that the following diagram commutes.
\[\begin{tikzcd}\ex L\rar{\phi_{0}}\dar&\et m{P_{v}}\dar
\\ \ex F\rar&\et {m-1}{P/v}\end{tikzcd}\]
In the case that $v$ is primitive, the above is a pullback diagram, and a family in this component of $\mathcal {BG}(\et mP)$ parametrized by $\ex F$ is equivalent to a choice of such a map $\ex F\longrightarrow \et {m-1}{P/v}$. The universal family over such a component is the projection $\et m{P_{v}}\longrightarrow \et {m-1}{P/v}$.
If $v$ is $k$ times a primitive integral-vector, the $v$--component of $\mathcal {BG}(\et mP)$ is  the quotient of $\et {m-1}{P/v}$ by the trivial $\mathbb Z_{k}$ action. In this case, the universal family over $\et {m-1}{P/v}$ has tropical part a polytope $P'_{v}$ with a given integral-affine map $P'_{v}\longrightarrow P_{v}$ that stretches a primitive vector to $v$, and is an affine isomorphism with determinant $k$. There is a corresponding degree-$k$ map $\et m{P_{v}'}\longrightarrow \et m{P_{v}}$ that gives a family 
\[\begin{tikzcd}\et m{P'_{v}}\dar\rar& \et m{P_{v}}
\\ \et {m-1}{P/v}\end{tikzcd}\]
with automorphism group $\mathbb Z_{k}$ (multiplying fibers by $k$th roots of unity.) Given any other family $\ex L$ in this component of $\mathcal {BG}(\et mP)$, there is a unique $\mathbb Z_{k}$--fold cover, $\ex L\times_{\et m{P_{v}}}\et m{P'_{v}}$ of $\ex L$ with a $\mathbb Z_{k}$--equivariant morphism to this family.  Conversely, given a $\mathbb Z_{k}$--fold cover $\hat{\ex F}$ of $\ex F$ and a $\mathbb Z_{k}$--invariant map $\psi\co  \hat {\ex F}\longrightarrow \et {m-1}{P/v}$, we can construct a family parametrized by $\ex F$  in this component of $\mathcal {BG}(\et mP)$ by taking the pullback of the above family by $\psi$, then taking the quotient by the $\mathbb Z_{k}$--action.  In summary, this $v$--component of $\mathcal {BG}(\et mP)$ is the quotient of $\et {m-1}{P/v}$ by the trivial $\mathbb Z_{k}$--action, and the pullback of the universal family to $\et{m-1}{P/v}$ is the family displayed above.

\

Below, we shall often restrict interest to components of $\mathcal G(\ex B)$ corresponding to vectors that span infinite rays in $\totb{\ex B}$. We shall use the notation $\gend(\ex B)$ for these components of $\mathcal G(\ex B)$, and $\End(\ex B)$ for the  quotient stack of $\gend(\ex B)$.  The following lemma establishes that $\End(\ex B)$ is complete if $\ex B$ is complete; the same fails to hold for $\mathcal {BG}(\ex B)$ in general.

\begin{lemma}If $P$ is complete, and $v$ spans an infinite ray in $P$, then $P/v$ is complete.  
\end{lemma}

\pf As $P$ is complete, it is a subset of $\mathbb R^{m}$ defined by some finite set of inequalities $\alpha_{i}\geq 0$ where $\alpha_{i}\co  \mathbb R^{m}\longrightarrow \mathbb R$ is integral-affine. $P_{v}$ consists of the union of all strata tangent to $v$, so $P_{v}$ is defined by the equations $\alpha_{i}\geq 0$ and $\alpha_{j}>0$ for all such $\alpha_{j}$ so that $v\alpha_{j}\neq 0$. As $v$ spans an infinite ray in $P$,  $v\alpha_{i}\geq 0$ for all $i$, therefore traveling in the direction of $v$, we can make all $\alpha_{j}$ so that $v\alpha_{j}\neq 0$ as large as we like. It follows that $P/v\subset \mathbb R^{n}/v$ is defined by the projection of the equations $\alpha_{i}\geq 0$  for $\alpha_{i}$ so that $v\alpha_{i}=0$. Therefore, $P/v$ is closed, as required.

\stop 

\begin{lemma}\label{BG functor}$\mathcal {BG}$ is a functor from the category of exploded manifolds to  stacks over  the category of exploded manifolds. 
\end{lemma}

\pf
As $\mathcal G$ is a functor to the category of Lie groupoids, this lemma follows from the functoriality  of passing from Lie groupoids to quotient stacks,  given in \cite{orbifoldstack}. Still, we should  check  this still works when  exploded manifolds are used instead of smooth manifolds. 

Given a map of exploded manifolds, $\psi\co  \ex A\longrightarrow \ex B$, and a family 
\[\begin{tikzcd}\ex L\times_{\ex F}\ex L\dar\rar{\phi_{1}}&\mathcal G_{1}\ex A\dar
\\\ex L\dar\rar{\phi_{0}}&\mathcal G_{0}\ex A
\\ \ex F\end{tikzcd}\]
in $\mathcal {BG}(\ex A)$, composing $\phi_{i}$ with $\mathcal G_{i}\psi$ may not give a family in $\mathcal {BG}(\ex B)$, because $\ex L\times_{\ex F}\ex L$ may not be the pullback of $\mathcal G_{1}\ex B$, or equivalently in this case, $(\mathcal G_{0}\psi)\circ\phi_{0}$  may not be complete restricted to each fiber of $\ex L\longrightarrow \ex F$ . There is, however, a unique extension of this bundle to a family in $\mathcal {BG}(\ex B)$ constructed by lengthening the fibers of $\ex L\longrightarrow \ex F$ until the map to $\mathcal G_{0}(\ex B)$ is fiberwise complete. This is equivalently given by replacing $\ex L$ by the quotient of $\ex L\times_{\mathcal G_{0}\ex B}\mathcal G_{1}(\ex B)$ by the diagonal action of $\mathcal G(\ex A)$ on $\ex L$ and $\mathcal G_{1}\ex B$; see \cite{orbifoldstack}, Remarks 3.24 and 3.30 and the beginning of section 4.  In our case, this means identifying $(l,c)\in \ex L\times_{\mathcal G_{0}(\ex B)}\mathcal G_{1}\ex B$ with $(c'*l,c'^{-1}c)$ for all $c'$ so that $c'*l$ is defined. The map to $\mathcal G_{0}\ex B$ sends $(l,c)$ to $c*(\mathcal G_{0}\psi\circ \phi_{0})(l)$. The map from the square of this bundle to $\mathcal G_{1}(\ex B)$ sends $(l,c_{1},c_{2}*l,c_{3})$ to $(c_{1}*(\mathcal G_{0}\psi\circ \phi_{0})(l), c_{1}^{-1}c_{2}c_{3})$.  All this really does in our case is lengthen the fibers of $\ex L$ until  $\ex L\times_{\ex F}\ex L$ is the pullback of $\mathcal G_{1}(\hat {\ex B})$. The fact that such an extension is unique implies that this construction is compatible with morphisms between families in $\mathcal {BG}(\ex A)$ so we obtain a functor $\mathcal {BG}(\ex A)\longrightarrow \mathcal {BG}(\ex B)$. Similarly, this uniqueness implies that this construction is compatible with composition of maps of exploded manifolds, so $\mathcal {BG}$ defines a functor from the category of exploded manifolds to the category of stacks, as required. 
\stop

\begin{remark}Some tropical information is discarded passing from $\mathcal G$ to $\mathcal {BG}$ so that it is not true in general that $\mathcal {BG}(\ex A\times _{C}\ex B)$ is equal to $\mathcal {BG}(\ex A)\times_{\mathcal {BG}(\ex C)}\mathcal {BG}(\ex B)$. An example is given by $\ex A=\et1{(0,\infty)}\subset \ex T=\ex C$, $\ex B=\et 1{(-\infty,0)}\subset \ex T$. This loss of tropical information is reflected in the need to lengthen the fibers of $\ex L$ in the proof of Lemma \ref{BG functor} above. 
\end{remark}

\begin{defn}Let $\End(\ex B)\subset \mathcal{BG}(\ex B)$ be the substack consisting of connected components of $\mathcal {BG}(\ex B)$ corresponding to vectors spanning infinite rays in $\totb{\ex B}$, and let $\gend(\ex B)\subset\mathcal G(\ex B)$ be the corresponding Lie groupoid. 

Given an integral-vector $w$ in  ${}^{\mathbb Z}T\ex B$ or the tangent space to $\totb{\ex B}$, let $\mathcal G^{w}_{0}(\ex B)\subset{}^{\mathbb Z} T\ex B$ be the closure of the stratum containing $w$, and let $\mathcal G^{w}(\ex B)$ be the pullback of $\mathcal G(\ex B)$ via the inclusion $\mathcal G^{w}_{0}(\ex B)\subset \mathcal G_{0}(\ex B)$. 
Given a tropical curve $\gamma$ in $\totb{\ex B}$, choose a primitive integral-vector  on each edge $e$ of the domain of $\gamma$, and let $w_{e}$ be its image in $T\totb{\ex B}$. Then define the following Lie groupoid
\[\gend^{[\gamma]}(\ex B):=\prod_{e}\mathcal G^{w_{e}}(\ex B)\] 
and let $\End^{[\gamma]}(\ex B)$ be the quotient stack of $\gend^{[\gamma]}(\ex B)$.
\end{defn}

Note that there is a canonical isomorphism  $\mathcal G^{w}(\ex B)\longrightarrow \mathcal G^{-w}(\ex B)$ sending $(v,c)$ to $(-v,c^{-1})$, so our construction of $\gend^{[\gamma]}(\ex B)$ does not depend on orientation of our  primitive integral-vector on each edge. For compatibility with $\gend(\ex B)$, we always choose the outgoing primitive integral-vector on each end of $\gamma$.

Each connected component of $\End(\ex B)$ or $\End^{[\gamma]}(\ex B)$ is the quotient of an exploded manifold by a trivial group action. Forgetting this group action gives an exploded manifold $\rend(\ex B)$ or $\rend^{[\gamma]}(\ex B)$ respectively. For example, \[\rend(\et mP)=\et mP\coprod_{v}\et m{P/v}\]
where the disjoint union is over all nonzero integral-vectors $v$ that span an infinite ray in $P$. On the other hand, 
\[\End(\et mP)=\et mP/\ex T\coprod_{v}\et m{P/v}/\mathbb Z_{\abs v}\]
where $\abs v$ is the positive integer so that $v/\abs v$ is a primitive integral-vector, and the action of $\mathbb Z_{\abs v}$ and $\ex T$ is the trivial action.

\

%

 For each end of $\gamma$, there is an outward-pointing integral-vector on the corresponding edge of $\gamma$, and a corresponding map $\gend^{[\gamma]}(\ex B)\longrightarrow \gend(\ex B)$. If $\gamma$ has $n$ ends, the product of these maps gives a map 
\[\tilde i^{[\gamma]}\co  \gend^{[\gamma]}(\ex B)\longrightarrow (\gend(\ex B))^{n}\]
which induces the following map on the level of quotient stacks.
\[i^{[\gamma]}\co  \End^{[\gamma]}(\ex B)\longrightarrow (\End(\ex B))^{n}\]
So long as the domain of $\gamma$ is not $\mathbb R$, this map $i^{[\gamma]}$ is a projection (onto the closure of some stratum) that forgets the components $\mathcal G^{w_{e}}(\ex B)$ for internal edges $e$. The tropical completion of this map at the image of $\gamma$ will feature in our gluing formula, but to reduce notational complexity, we shall refer to this tropical completion again as $i^{[\gamma]}$.

The definition of $\gend^{[\gamma]}(\ex B)$ also makes sense if $\gamma$ is a cut tropical curve. Cut $\gamma$ at some chosen point on all edges, and let $\gamma^{+}_{v}$ be the connected component containing $v$. As each edge in $\gamma_{v}^{+}$ is part of an edge in $\gamma$, there is a canonical map 
\[\tilde \Delta\co  \gend^{[\gamma]}(\ex B)\longrightarrow \prod_{v}\gend^{[\gamma_{v}^{+}]}(\ex B)\ .\]

Now suppose that $\ex B$ is basic, and let $\bar{\ex B}_{v}$ be the closure of the stratum of $\ex B$ containing $v$. As $\ex B$ is basic, there is an inclusion  $\bar{\ex B}_{v}\subset\ex B\tc v$, where $\ex B\tc v$  is the tropical completion of $\ex B$ at $v$ discussed in \cite{vfc}. The  cut tropical curve $\gamma_{v}^{+}$ is contained in $\bar{\ex B}_{v}$, and $\gend^{[\gamma_{v}^{+}]}(\ex B)$ is a sub-groupoid of $\mathcal G(\bar{\ex B}_{v})^{n}$, so the inclusion $\bar{\ex B}_{v}\subset \ex B\tc v$ induces an inclusion $\gend^{[\gamma_{v}^{+}]}\longrightarrow \mathcal G(\ex B\tc v)^{n}$.  There is a tropical curve $\gamma_{v}$ in $\totb{\ex B\tc v}$ with a single vertex $v$ and all edges semi-infinite rays in the directions of the edges leaving $v$. The above inclusion has image contained in $\gend^{[\gamma_{v}]}(\ex B\tc v)$, so our map $\tilde \Delta$ above induces a canonical map
\[\tilde \Delta\co  \gend^{[\gamma]}(\ex B)\longrightarrow \prod_{v}\gend^{[\gamma_{v}]}(\ex B\tc v)\]
which induces the following map on the level of quotient stacks.
\[\Delta\co  \End^{[\gamma]}(\ex B)\longrightarrow \prod_{v}\End^{[\gamma_{v}]}(\ex B\tc v)\]
Again, it is the tropical completion of $\Delta$ at the image of $\gamma$ that will feature in our gluing formula, but to reduce notational complexity, we shall refer to this tropical completion simply as $\Delta$.

We now describe our evaluation maps to $\End(\ex B)$. Given a connected family of curves $\hat f$ in $\Mu_{\cdot}\ex B$ and an end-label, $e$, let $\tilde {\ex F}^{e}(\hat f)$ be the connected component of $\gend \ex C(\hat f)$ corresponding to the outgoing primitive integral-vector on the end $e$. Note that the quotient stack of $\tilde {\ex F}^{e}(\hat f)$ is canonically isomorphic to $\ex F(\hat f)$, so the restriction of $\End\hat f\co  \End(\ex C(\hat f))\longrightarrow \End\ex B$ to this component defines a natural evaluation map $\ex F(\hat f)\longrightarrow \End\ex B$. As $\End$ is a functor,   we get a corresponding evaluation map $\Mu_{\cdot}(\ex B)\longrightarrow \End \ex B$ for each end $e$.
The product of  these evaluation maps for all labeled ends gives an evaluation map  $\ev$, and composing this evaluation map with the map forgetting stack structure gives a further evaluation map,  $ev$.
\[\begin{tikzcd}\Mu_{\cdot}(\ex B)\ar[bend left]{rr}{ev}\rar{\ev}&\coprod_{n} (\End \ex B)^{n}\rar& \coprod_{n}(\rend \ex B)^{n}\end{tikzcd}\]

Let us lift this evaluation map $\ev$  to give some kind of evaluation map to the corresponding groupoid $\coprod_{n}(\gend \ex B)^{n}$.
Given a family of curves $\hat f$ in $\Mu_{\cdot}\ex B$, let $\tilde {\ex F}(\hat f)$ be the fiber-product of $\tilde{\ex F}^{e}(\hat f)$ over $\ex F(\hat f)$ for all ends $e$. Then $\tilde {\ex F}$ defines a functor from $\Mu_{\cdot}$ to the category of (exploded) Lie groupoids. The quotient stack of the groupoid $\tilde{\ex F}(\hat f)$ is equal to $\ex F(\hat f)$, and $\ev\co  \ex F(\hat f)\longrightarrow \coprod_{n}(\End\ex B)^{n}$ lifts to a map of Lie groupoids induced from $\mathcal G\hat f$.
\[\tilde{\ev}\co  \tilde {\ex F}(\hat f)\longrightarrow \coprod_{n}(\gend\ex B)^{n}\] 
Clearly, $\tilde{\ev}$ is a natural transformation from the functor $\tilde{\ex  F}$. We can also consider $\tilde {\ev}$ as a map of groupoids internal to the category of stacks over exploded manifolds.  Let $\Mu_{\cdot^{+}}$ be the moduli stack of curves in $\Mu_{\cdot}$ along with a choice of point in every external edge. We can consider $\tilde {\ev}$ as given by the maps $\tilde {\ev}_{i}$ of stacks in the following commutative diagram
\[\begin{tikzcd}\ex F(\hat f^{++})=\tilde{\ex F}_{1}(\hat f)\dar[shift left]\dar[shift right]\rar&\Mu_{\cdot^{++}}\dar[shift left]\dar[shift right]\rar{\tilde{\ev}_{1}} & \coprod_{n}(\gend_{1}\ex B)^{n}\dar[shift left]\dar[shift right]
\\\ex F(\hat f^{+})=\tilde{\ex F}_{0}(\hat f)\dar\rar &\Mu_{\cdot^{+}}\dar\rar{\tilde{\ev}_{0}} &\coprod_{n}(\gend_{0}\ex B)^{n}\dar
\\ \ex F(\hat f)\rar& \Mu_{\cdot}\rar{\ev}\ar[bend right]{rr}{ev}&\coprod_{n}(\End\ex B)^{n}\rar& \coprod_{n}(\rend \ex B)^{n}\end{tikzcd}\]
where $\Mu_{\cdot^{++}}:=\Mu_{\cdot^{+}}\times_{\Mu_{\cdot}}\Mu_{\cdot^{+}}$.
The map $\tilde{\ev}_{0}$ is given by evaluation of $T\hat f$ at the outgoing primitive integral-vector at each of the points chosen on external edges.  A curve in $\Mu_{\cdot^{++}}$ has two chosen points in each external edge. In a family $\hat f$, this corresponds to two sections $\ex F(\hat f)\longrightarrow \ex C(\hat f)$ for each external edge $e$; the second section is $c_{e}$ times the first section for some exploded function $c_{e}\co  \ex F(\hat f)\longrightarrow \ex T$. The map $\tilde {\ev}_{1}$ is the lift of $\tilde {\ev}_{0}$ determined by these functions $c_{e}$ when we consider $\gend_{1}(\ex B)\subset \gend_{0}(\ex B)\times \ex T$.  If we consider $\Mu_{\cdot^{++}}\rightrightarrows \Mu_{\cdot^{+}}$ a groupoid $\tilde{\Mu_{\cdot}}$ (internal to the category of stacks over exploded manifolds), the above maps $\tilde {\ev}_{i}$ can be represented as a map
\[\tilde{\ev}\co  \tilde {\Mu_{\cdot}}\longrightarrow \coprod_{n}(\gend\ex B)^{n}\]

Given a tropical curve $\gamma$ in $\totb{\ex B}$, we can describe  similar evaluation maps, $\ev^{[\gamma]}$, $ev^{[\gamma]}$ and $\tilde{\ev}^{[\gamma]}$, with target $\End^{[\gamma]}\ex B$, $\rend^{[\gamma]}\ex B$, and $\gend^{[\gamma]}\ex B$. In this case, we need to use the stack $\Mu_{[\gamma]}$ from Definition \ref{decorated}.  Recall from Definition \ref{stack+} and the discussion following it, that $\Mu_{[\gamma]^{+}}$ is the stack of curves in $\Mu_{[\gamma]}$ along with a choice of point in each edge decorated by an edge of $\gamma$, and that a family $\hat f$ in $\Mu_{[\gamma]}$ has a canonical lift to a family $\hat f^{+}$ in $\Mu_{[\gamma]}$ parmametrized by the fiber-product over $\ex F(\hat f)$ of the strata $\ex C_{e}\subset \ex C(\hat f)$ decorated by edges of $\gamma$. Define $\tilde{\ex F}_{0}(\hat f):=\ex F(\hat f^{+})$ and let $\tilde {\ex F}(\hat f)$ be the banal groupoid of the bundle $\tilde{\ex F}_{0}(\hat f)\longrightarrow \ex F(\hat f)$, so $\tilde{\ex F}_{1}=\tilde {\ex F}_{0}\times _{\ex F}\tilde {\ex F}_{0}$. If $\gamma$ has $n$ edges, choosing an orientation for each edge of $\gamma$ identifies $\tilde {\ex F}(\hat f)$ as a sub-groupoid of $(\mathcal G \ex C(\hat f))^{n}$. Then $\mathcal G{\hat f}$ applied to this sub-groupoid  has image in $\gend^{[\gamma]}\ex B\subset (\mathcal G \ex B)^{n}$, and defines our map $\tilde{\ev}^{[\gamma]}$.
\[\tilde{\ev}^{[\gamma]}\co  \tilde{\ex F}(\hat f)\longrightarrow \gend^{[\gamma]}\ex B\]
The map $\tilde{\ev}_{0}$ we have seen before. After noting that $\gend^{[\gamma]}_{0}(\ex B)=\prod_{e}\bar{\ex B}_{e}$, $\tilde {\ev}_{0}$ is equal to $\prod_{e} ev_{e}\co  \ex F(\hat f^{+})\longrightarrow \prod_{e}\bar{\ex B}_{e}$, where $ev_{e}$ first appears in Definition \ref{stack+}, and stars in Theorem \ref{stacks fiber-product}.
As the quotient groupoid of $\tilde {\ex F}(\hat f)$ is $\ex F(\hat f)$, the corresponding map on quotient groupoids defines an evaluation map $\ev^{[\gamma]}$.
\[\ev^{[\gamma]}\co  \ex F(\hat f)\longrightarrow \End^{[\gamma]}\ex B\]
As with $\tilde {\ev}$, we may consider $\tilde {\ev}^{[\gamma]}_{i}$ as giving maps from moduli stacks as in the following commutative diagram. 
\[\begin{tikzcd}\tilde{\ex F}_{1}(\hat f)\dar[shift left]\dar[shift right]\rar&\Mu_{[\gamma]^{++}}\dar[shift left]\dar[shift right]\rar{\tilde{\ev}^{[\gamma]}_{1}} & \gend^{[\gamma]}_{1}\ex B\dar[shift left]\dar[shift right]
\\\tilde{\ex F}_{0}(\hat f)\dar\rar &\Mu_{[\gamma]^{+}}\dar\rar{\tilde{\ev}^{[\gamma]}_{0}} &\gend^{[\gamma]}_{0}\ex B\dar
\\ \ex F(\hat f)\rar& \Mu_{[\gamma]}\rar{\ev^{[\gamma]}}\ar[bend right]{rr}{ev^{[\gamma]}}&\End^{[\gamma]}\ex B\rar& \rend^{[\gamma]}\ex B\end{tikzcd}\]
In the above, $\Mu_{[\gamma]^{++}}:=\Mu_{[\gamma]^{+}}\times_{\Mu_{[\gamma]}}\Mu_{[\gamma]^{+}}$ is the moduli stack of curves in $\Mu_{[\gamma]}$ with two chosen points on each edge labeled by an edge of $\gamma$. We may use $\tilde{\ev}^{[\gamma]}\co  \tilde{\Mu_{[\gamma]}}\longrightarrow \gend^{[\gamma]}(\ex B)$ as a shorthand for the top-right square of the above diagram.

\

All curves in $\Mu_{[\gamma_{v}^{+}]}(\ex B)$ have image contained in $\bar {\ex B}_{v}$, so the inclusion $\bar{\ex B}_{v}\subset\ex B\tc v$ induces an inclusion $\Mu_{[\gamma_{v}^{+}]}(\ex B)\subset \Mu_{[\gamma_{v}^{+}]}(\ex B\tc v)$.
Theorem \ref{stacks fiber-product} along with the observation that $\gamma_{v}^{+}$ has  no internal edges and has all external edges infinitely extendable,  implies that the cutting map applied to $\Mu_{[\gamma_{v}]^{+}}(\ex B\tc v)$ is an isomorphism, so 
\[\Mu_{[\gamma_{v}^{+}]}(\ex B\tc v)=\Mu_{[\gamma_{v}]^{+}}(\ex B\tc v)\ .\]
 Define $\cut_{0}$ to be the following composition.
\[\begin{tikzcd}\Mu_{[\gamma]^{+}}(\ex B)\rar{\cut}\ar[bend left]{rr}{\cut_{0}}&\prod_{v}\Mu_{[\gamma_{v}^{+}]}(\ex B)\rar[hook]&\prod_{v}\Mu_{[\gamma_{v}]^{+}}(\ex B\tc v)\end{tikzcd}\]

Our map $\cut_{0}$ determines maps, $\cut_{1}$ and $\cut$, so that  the following diagram commutes.
\begin{equation}\label{cutcd}\begin{tikzcd}\Mu_{[\gamma]^{++}}(\ex B)\dar{\cut _{1}}\rar[shift left]{s}\rar[shift right,swap]{t} &\Mu_{[\gamma]^{+}}(\ex B)\rar{\pi}\dar{\cut_{0} }&\Mu_{[\gamma]}(\ex B)\dar{\cut }
\\ \prod_{v}\Mu_{[\gamma_{v}]^{++}}(\ex B\tc v)\rar[shift left]{s}\rar[shift right, swap]{t} &\prod_{v}\Mu_{[\gamma_{v}]^{+}}(\ex B\tc v)\rar{\pi}&\prod_{v}\Mu_{[\gamma_{v}]}(\ex B\tc v)\end{tikzcd} \end{equation}
So long as $\gamma$ has internal edges, $(\cut\hat f)^{+}$ is larger than $\cut_{0}(\hat f^{+})$, because $\cut\hat f$ has two edges for each internal edge of $\gamma$. The universal property of $(\cut\hat f)^{+}$ and $(\cut\hat f)^{++}$ give maps $\cut_{0}(\hat f^{+})\longrightarrow (\cut_{0} \hat f)^{+}$ and $\cut_{1}(\hat f^{++})\longrightarrow (\cut \hat f)^{++}$ that define a map of groupoids $\tilde{\ex F}(\hat f)\longrightarrow \tilde{\ex F}(\cut\hat f)$. 

\

Use the notation $\Mu_{[\gamma_{v}]}(\ex B)\subset \Mu_{[\gamma_{v}]}(\ex B\tc v)$ for the substack consisting of curves whose tropical part has all vertices in $\totb{\bar{\ex B}_{v}}\subset \totb{\ex B\tc v}$. This is the image of $\Mu_{[\gamma_{v}^{+}]}(\ex B)\subset \Mu_{[\gamma_{v}]^{+}}(\ex B\tc v)$ under the projection forgetting the extra points on edges. From Theorem \ref{stacks fiber-product},  one could expect that the following diagram
\[\begin{tikzcd}\Mu_{[\gamma]}(\ex B)\rar{\ev^{[\gamma]}}\dar{\cut }&\End^{[\gamma]}\ex B\dar{\Delta}
\\ \prod_{v}\Mu_{[\gamma_{v}]}(\ex B)\rar{\prod_{v}\ev^{[\gamma_{v}]}}&\prod_{v}\End^{[\gamma_{v}]}\ex B\tc v
\end{tikzcd}\]
is close to being a fiber-product diagram. It fails to be a fiber-product diagram because, in taking the quotient, we have thrown away tropical information encoding the requirement that the length of an internal edge of a tropical curve is positive. We shall deal with this issue using tropical completion, but shall do so on the level of Kuranishi structures instead of using $\Mu_{[\gamma]}$.

\begin{prop}\label{cut*} Given any family of curves $\hat f$ in $\prod_{v}\Mu_{[\gamma_{v}]}(\ex B)$ so that $\prod_{v}ev^{[\gamma]}$ is transverse to $\Delta$, there exists a family of curves,  $\cut^{*}  (\hat f)$, in $\Mu_{[\gamma]}(\ex B)$ with a map 
\[\cut (\cut^{*}  (\hat f))\longrightarrow \hat f\] 
satisfying the following universal property: Given any family of curves $\hat h$ in $\Mu_{[\gamma]}(\ex B)$ and a map $\psi\co  \cut (\hat h)\longrightarrow \hat f$, there exists a unique map $\cut^{*}  (\psi)\co  \hat h\longrightarrow \cut^{*}  (\hat f)$ so that the following diagram commutes.
\begin{equation}\label{universal glue}\begin{tikzcd} & \cut (\cut^{*}  (\hat f))\dar
\\ \cut (\hat h)\ar{ur}{\cut (\cut^{*}  (\psi))}\rar{\psi}& \hat f
\end{tikzcd}\end{equation}
Moreover,  the following diagram
\begin{equation}\label{gpd fp}\begin{tikzcd}\tilde{\ex F}(\cut^{*}  (\hat f))\dar \rar{\tilde{\ev}^{[\gamma]}}&\gend^{[\gamma]}(\ex B)\dar{\tilde\Delta}
\\ \tilde{\ex F}(\hat f)\rar{\prod_{v}\tilde{\ev}^{[\gamma_{v}]}}&\prod_{v}\gend^{[\gamma]}(\ex B\tc v)\end{tikzcd}\end{equation}
is a fiber-product diagram of groupoids.

\end{prop}

\pf

Theorem \ref{stacks fiber-product} allows us to construct $\cut_{0}^{*}(\hat f^{+})$  so that the following is a fiber-product diagram.
\[\begin{tikzcd}\ex F(\cut_{0}^{*}\hat f^{+})\dar\rar{\tilde{\ev}^{[\gamma]}_{0}}&\gend^{[\gamma]}_{0}(\ex B)\dar{\tilde\Delta_{0}}
\\ \ex F(\hat f^{+})\rar{\prod_{v}\tilde{\ev}_{0}^{[\gamma_{v}]}}& \prod_{v}\gend^{[\gamma_{v}]}(\ex B\tc v)\end{tikzcd}\]
In particular, define $\cut_{0}^{*}(\hat f^{+})$ as the construction, $\cut^{*}$, from Theorem \ref{stacks fiber-product} applied to the fiber-product of $\hat f^{+}$ with $\gend^{[\gamma]}_{0}(\ex B)$: taking this fiber-product gives a family of curves in $\Mu_{[\gamma_{v}]^{+}}(\ex B\tc v)$ which, when cut, are contained entirely inside $\bar{\ex B}_{v}\subset \ex B\tc v$, so correspond to a family of curves in $\Mu_{[\gamma_{v}^{+}]}(\ex B)$ to which Theorem \ref{stacks fiber-product} applies. The fiber-product exists because $\prod_{v}ev^{[\gamma_{v}]}$ being transverse to $\Delta$ implies that $\prod_{v}\tilde{\ev}^{[\gamma_{v}]}_{0}$ is transverse to $\tilde \Delta_{0}$. 

 The universal property from Theorem \ref{stacks fiber-product} implies the corresponding universal property for $\cut^{*}_{0}\hat f^{+}$: given any family of curves $\hat h$ in $\Mu_{[\gamma]^{+}}(\ex B)$, any map $\psi\co  \cut_{0}\hat h\longrightarrow \hat f^{+}$ has a unique pullback, $\cut_{0}^{*}\psi\co  \hat h\longrightarrow \cut_{0}^{*}(\hat f^{+})$ so that the following diagram commutes.
 \[\begin{tikzcd} & \cut_{0} (\cut_{0}^{*}  (\hat f^{+}))\dar
\\ \cut (\hat h)\ar{ur}{\cut_{0} (\cut_{0}^{*}  (\psi))}\rar{\psi}& \hat f^{+}\end{tikzcd}\]

We shall construct $\cut^{*}\hat f$ so that $\cut_{0}^{*}(\hat f^{+})=(\cut\hat f)^{+}$ (or rather these two families are canonically isomorphic, because of the universal property of pullbacks and the commutative diagram (\ref{cutcd})). Using the notation of diagram (\ref{cutcd}), $\pi\circ\cut_{0}\circ s=\pi\circ\cut_{0}\circ t=\pi \circ s\circ\cut_{1}$, so there are canonical isomorphisms (written as equalities) as follows:
\[s^{*}\cut_{0}^{*}\hat f^{+}=t^{*}\cut_{0}^{*}\hat f^{+}:=\cut_{1}^{*}(\hat f^{++})\]
The natural maps
\[ss^{*}\cut_{0}^{*}\hat f^{+}\longrightarrow \cut^{*}_{0}\hat f^{+} \]
\[tt^{*}\cut_{0}^{*}\hat f^{+}\longrightarrow \cut^{*}_{0}\hat f^{+}\]
project to give the following groupid internal to $\Mu_{[\gamma]}(\ex B)$.
\[\pi s\cut_{1}^{*}(\hat f^{++})\rightrightarrows\pi \cut_{0}^{*}(\hat f^{+})\]
Define $\cut^{*}\hat f$ to be the quotient of this groupoid. (The groupoid action is free, so  $\ex F(\cut^{*}\hat f)$ is an exploded manifold instead of a stack.) The map $\pi\cut_{0}^{*}\hat f^{+}\longrightarrow \cut^{*}\hat f$ lifts uniquely to a map $\cut_{0}^{*}\hat f^{+}\longrightarrow(\cut ^{*}\hat f)^{+}$ with inverse the unique lift of the map $\pi\circ \cut_{0}(\cut^{*}f)^{+}\longrightarrow f$ using the universal property of $\cut_{0}^{*}$ and $\pi^{*}$. So $(\cut^{*}\hat f)^{+}$ is indeed canonically isomorphic to $\cut_{0}^{*}\hat f^{+}$.

The universal property for $\cut^{*}\hat f$ is implied by the universal property for $\cut_{0}^{*}\hat f^{+}$. In particular, given any $\hat h$ in $\Mu_{[\gamma]^{+}}(\ex B)$, any map $\psi\co  \cut\pi\hat h\longrightarrow \hat f$ has a unique lift to $\psi^{+}\co  \cut_{0}\hat h\longrightarrow \hat f^{+}$ and $\cut_{0}^{*}\psi^{+}\co  \hat h\longrightarrow \cut_{0}^{*}\hat f^{+}=(\cut^{*}\hat f)^{+}$. Projecting this map using $\pi$ gives the required unique lift $\cut^{*}\psi\co  \pi\hat h\longrightarrow \cut^{*}\hat f$. As all families in $\Mu_{[\gamma]}(\ex B)$ are locally in the image of $\pi$, the universal property of $\cut^{*}$ holds for all families. 

It remains to show that diagram (\ref{gpd fp}) is a fiber-product diagram. We have already identified $\ex F(\cut_{0}^{*}\hat f^{+})=\tilde {\ex F}_{0}(\cut^{*}\hat f)$ as the appropriate fiber-product, so it remains to check that 
\begin{equation}\label{gpd fp1}\begin{tikzcd}\tilde{\ex F}_{1}(\cut^{*} \hat f)\dar \rar{\tilde{\ev}_{1}^{[\gamma]}}&\gend_{1}^{[\gamma]}(\ex B)\dar{\tilde\Delta}
\\ \tilde{\ex F}_{1}(\hat f)\rar{\prod_{v}\tilde{\ev}_{1}^{[\gamma_{v}]}}&\prod_{v}\gend^{[\gamma]}_{1}(\ex B\tc v)\end{tikzcd}\end{equation}
is a fiber-product diagram. A map $x\co  \ex A\longrightarrow \tilde{\ex F}_{1}(\cut^{*}\hat f)$ is equivalent to a map $x_{0}\co  \ex A\longrightarrow \tilde{\ex F}_{0}(\cut^{*}\hat f)$ and for each edge $e$ of $\gamma$, a map $c_{e}\co  \ex A\longrightarrow \ex T$ so that the action of $c_{e}$ on $x_{0}^{*}(\cut^{*}\hat f)^{+}$ is well-defined.  We already know that $x_{0}$ is equivalent to the corresponding maps to $\gend_{0}(\ex B)$ and $\tilde{\ex F}_{0}(\hat f)$ compatible with the fiber-product. Promoting these to maps to compatible maps to $\gend_{1}(\ex B)$ and $\tilde{\ex F}_{1}(\hat f)$ is equivalent to choosing functions $c_{e}\co  \ex A\longrightarrow \ex T$ for every edge of $\gamma$ and $\coprod_{v}\gamma_{v}$ so that
\begin{enumerate}
\item \label{p1}$c_{e_{1}}=c^{-1}_{e_{2}}=c_{e}$ if $e$ is an internal edge of $\gamma$ that breaks into $e_{1}$ and $e_{2}$.
\item\label{p2} $c_{e}=c_{e'}$ if $e'$ is the edge of $\coprod_{v}\gamma_{v}$ corresponding to an end, $e$, of $\gamma$. 
\item\label{p3} For edges $e$ of $\coprod_{v}\gamma_{v}$, the action of $c_{e}$ on the pullback of $\hat f^{+}$ is well-defined.
\item \label{p4}For edges $e$ of $\gamma$, the action of $c_{e}$ on the image of $\ex A$ in $\gend^{[\gamma]}_{0}{\ex B}$ is well-defined.

\end{enumerate}
These conditions on our functions are equivalent to a choice of $c_{e}$ for each edge of $\gamma$ so that the action on $x_{0}^{*}(\cut^{*}\hat f)^{+}$ is well-defined:  Conditions (\ref{p1}) and (\ref{p2}) ensure that we only need $c_{e}$ for edges of $\gamma$; condition (\ref{p3}) (along with the identifications from the previous conditions) is equivalent to the action of $c_{e}$ being well-defined on 
$x_{0}^{*}(\cut^{*}\hat f)^{+}$, and condition (\ref{p4}) follows from the previous conditions and the fact that the maps to $\gend_{0}(\ex B)$ and $\tilde{\ex F}_{0}(\hat f)$ are compatible. Therefore, diagram (\ref{gpd fp1}) is a fiber-product diagram, and the proof that diagram (\ref{gpd fp}) is a fiber-product diagram is complete.

\stop

\section{Simple gluing formula}\label{simple gluing formula}

Proposition \ref{cut*} implies that the following is almost a fiber-product diagram
\[\begin{tikzcd}\Msw_{[\gamma]}(\ex B)\dar{\cut}\rar{\ev^{[\gamma]}}&\End^{[\gamma]}(\ex B)\dar{\Delta}
\\ \prod_{v}\Msw_{[\gamma_{v}]}(\ex B\tc v)\rar{\prod_{v}\ev^{[\gamma_{v}]}}&\prod_{v}\End^{[\gamma_{v}]}(\ex B\tc v)\end{tikzcd}\]
so we could hope that a version of the usual relationship for pushforwards of differential forms holds. It does hold after suitable tropical completion, and results in our gluing formula. Of course, $\ev^{[\gamma]}_{!}$, as constructed in \cite{vfc}, involves pushing forward from the virtual fundamental class, which is  the intersection with $0$ of a section $\dbar$
 of a sheaf $\mathcal Y$ over $\Msw$. To prove our gluing formula, we must first verify that our diagram above is compatible with $\dbar$ and $\mathcal Y$.

 We now describe a natural identification of $\mathcal Y$ over $\Msw_{[\gamma]}(\ex B)$ with $\cut^{*}\mathcal Y$. Using the cutting map,
a family of curves $\hat f$ is sent to several families of curves,  $\hat f_{v}$, all parametrized by $\ex F(\hat f)$. There is an identification, of sheaves of $\C\infty1(\ex F(\hat f))$--modules, in the following form.
\[\mathcal Y(\hat f)\cong \oplus_{v}\mathcal Y(\hat f_{v})\]

Using the cutting map, the strata, $\ex C(\hat f)_{v}$,  of $\ex C(\hat f)$ labelled by $v$ have a natural inclusion into $\ex C(\hat f_{v})$, compatible with the inclusion,  $\bar{\ex B}_{v}\subset \ex B\tc v$, of the closure of the stratum of $\ex B$ containing $v$ into $\ex B\tc v$.
\begin{equation}\label{cut inclusions}\begin{tikzcd}\ex C(\hat f)_{v}\dar[hook]\rar{\hat f}&\bar{\ex B}_{v}\dar[hook]
\\ \ex C(\hat f_{v})\rar{\hat f_{v}}&\ex B\tc v\end{tikzcd}\end{equation}
The sheaf $\mathcal Y(\hat f_{v})$ consists of $\C\infty1$ sections of $T^{*}_{vert}\ex C(\hat f_{v})\otimes \hat f_{v}^{*}T\ex B\tc v$ that vanish on the edges of curves. Our inclusions above identify such sections with sections of $T^{*}_{vert}\ex C(\hat f)\otimes \hat f^{*} T\ex B$ over $\ex C(\hat f)_{v}$. We can extend these sections to be $0$ everywhere else on $\ex C(\hat f)$; these extended sections are $\C\infty1$ because  our sections vanish on edges.  This identification defines an inclusion $\mathcal Y(\hat f_{v})\subset \mathcal Y(\hat f)$. 

The resulting map $\oplus_{v}\mathcal Y(\hat f_{v})\longrightarrow \mathcal Y(\hat f)$ has an obvious inverse. Any section of $\mathcal Y(\hat f)$ restricted to $\ex C(\hat f)_{v}$, then extended as $0$ elsewhere on $\ex C(\hat f_{v})$ defines a section of $\mathcal Y(\hat f_{v})$. This construction for all $v$ provides our required inverse map,   $\mathcal Y(\hat f)\longrightarrow \oplus_{v}\mathcal Y(\hat f_{v})$. These identifications are clearly compatible with the $\C\infty 1(\ex F(\hat f))$--module structure and the section $\dbar$,  and are natural in the sense that given any map $\hat f\longrightarrow \hat g$ in $\Msw_{[\gamma]}(\ex B)$, the corresponding diagram commutes.
\[\begin{tikzcd}\mathcal Y(\hat f)\dar &\lar \mathcal Y(\hat g)\dar
\\ \oplus_{v}\mathcal Y(\hat f)&\lar \oplus_{v}\mathcal Y(\hat g_{v})\end{tikzcd}\]
This completes the description of our natural identification \[\mathcal Y(\Msw_{[\gamma]}(\ex B))=\cut^{*}(\mathcal Y(\prod_{v}\Msw_{[\gamma_{v}]}(\ex B\tc v)))\ .\] 

The virtual fundamental class of the moduli stack of holomorphic curves, constructed in \cite{vfc}, uses an embedded Kuranishi structure, constructed in \cite{evc}. Each Kuranishi chart on $\Msw_{\bullet}$ is an open substack $\mathcal U\subset\Msw_{\bullet}$, and on $\mathcal U$ a finite-dimensional obstruction bundle $V\subset \mathcal Y$, and a family $\hat f$ with automorphism group $G$ so that $\hat f/G$ represents $\dbar^{-1}(V)\subset \mathcal U$.

\begin{remark} Proposition \ref{cut*} implies  that if $\hat f/G$ represents the substack $\dbar^{-1}(V)$, then $(\cut^{*}\hat f)/G$ represents the substack $\dbar^{-1}(\cut^{*}V)$.  
 \end{remark}

 The pullback of a Kuranishi chart $(\mathcal U,V,\hat f/G)$ on $\prod_{v}\Msw_{[\gamma_{v}]}(\ex B_{v})$ is
 \[ \cut^{*}(\mathcal U,V,\hat f/G):=(\cut^{-1}\mathcal U,\cut^{*}V,(\cut^{*}\hat f)/G)\] however we need to check that our pulled-back Kuranishi chart satisfies the technical conditions of Definition 2.27 of \cite{evc}.
To be used for a Kuranishi chart as defined in \cite{evc}, an obstruction bundle must be simply-generated in the sense of Definition 2.24 of \cite{evc}.

\begin{lemma}\label{sg pullback} Let $f\in \Msw_{[\gamma]}(\ex B)$ be a holomorphic curve with $\cut f=\prod_{v} f_{v}$. Given any choice of simply-generated obstruction bundles $V_{v}$, defined on a open neighborhoods of $ f_{v}$ in $\Msw_{[\gamma_{v}]}(\ex B\tc v)$, the pullback, $\cut^{*}(\oplus V_{v})$,  is simply-generated on an open neighorhood of $f$ in $\Msw_{[\gamma]}(\ex B)$.
\end{lemma}

\pf

Lemma 6.4 from \cite{evc} implies that,  on a small enough neighborhood, $\mathcal U_{v}$ of $f_{v}$,   $V_{v}$ is pulled back (Definition 2.23 of \cite{evc}) from some family of curves $\hat g_{v}$ in $\mathcal U_{v}$ with a group $G_{v}$ of automorphisms using a map 
\[\begin{tikzcd}\mathcal U_{v}^{+1}\rar\dar&\ex C(\hat g_{v})/G_{v}\dar
\\\mathcal U_{v}\rar&\ex F(\hat g_{v})/G_{v}
\end{tikzcd}\]
that is a holomorphic isomorphism restricted to each fiber.  Consider   $\prod_{v}\hat g_{v}$ as a family of (possibly unstable) holomorphic curves, mapping to a point, and decorated by $\gamma_{v}$. Let $\hat g$ be $\cut^{*}(\prod_{v}\hat g_{v})$ in $\Mpsw_{[\gamma]}(pt)$. As specified by Proposition \ref{cut*}, this $\hat g$ comes with a natural morphism $\cut \hat g\longrightarrow \prod_{v}\hat g_{v}$, which satisfies the property that given any family of  curves $\hat f$ in $\Mpsw_{[\gamma]}(pt)$ with a map $\psi\co  \cut (\hat f)\longrightarrow \prod_{v}\hat g_{v}$, there exists a unique map $\psi'\co  \hat f\longrightarrow \hat g$ so that the following diagram commutes:
\[\begin{tikzcd}\cut (\hat f)\ar[bend left]{rr}{\psi}\rar{\cut (\psi')} & \cut \hat g\rar& \prod_{v}g_{v}\end{tikzcd}\] 
This universal property implies that the fiberwise holomorphic map 
\[\begin{tikzcd} \prod_{v}\mathcal U_{v}^{+1}\rar\dar&\prod_{v}\ex C(\hat g_{v})/G_{v}\dar
\\  \prod_{v}\mathcal U_{v}\rar&\prod_{v}\ex F(\hat g_{v})/G_{v}
\end{tikzcd}\]
lifts canonically to a fiberwise holomorphic map 
\begin{equation}\label{fhm}\begin{tikzcd}\mathcal U^{+1}\rar\dar&\ex C(\hat g)/\prod_{v}G_{v}\dar
\\\mathcal U\rar&\ex F(\hat g)/\prod_{v}G_{v}
\end{tikzcd}\end{equation}
where $\mathcal U$ is the inverse image of $\prod_{v}\mathcal U_{v}$ under the cutting map. We shall now show that using the above map (\ref{fhm}), $\cut^{*}(\oplus_{v}V_{v})$ is the pullback of a sheaf of sections in the sense of Definition 2.23 of \cite{evc}, and in particular $\cut^{*}(\oplus_{v} V_{v})$ is simply-generated in the sense of Definition 2.24.  As stipulated by these definitions, $V_{v}$ is constructed using the pullback of some sections $\sigma$ of $\Gamma^{0,1}(T^{*}_{vert}\ex C(\hat g_{v})\otimes T(\ex B\tc v))$ vanishing on edges of $\ex C(\hat g_{v})$. Using the inclusions (\ref{cut inclusions}), and the map $\ex C(\cut \hat g)\longrightarrow \ex C(\prod_{v}\hat g_{v})$, any such section $\sigma$ can be pulled back, then  extended to a section $\sigma'$ of $\Gamma^{0,1}(T^{*}_{vert}\ex C(\hat g))\otimes T\ex B)$, vanishing outside the strata of $\ex C(\hat g)$ labeled by $v$.  These sections $\sigma'$, when pulled back\footnote{Definition 2.23 of \cite{evc}.} using our map (\ref{fhm}) generate $\cut^{*}(\oplus_{v}V)$  in the sense of 2.24 of \cite{evc}.

This completes the proof that the pullback of $\oplus_{v}V_{v}$ is simply-generated on $\mathcal U$.

\stop

 \
 
\

\begin{prop}\label{Kuranishi fiber-product} Given any choice of $ev^{[\gamma_{v}]}$--submersive  embedded Kuranishi structure on $\Msw_{[\gamma_{v}]}(\ex B\tc v)$, there exists an embedded Kuranishi structure on $\prod_{v}\Msw_{[\gamma_{v}]}(\ex B\tc v)$   so that every Kuranishi chart has an extension that is a product of charts from the chosen Kuranishi structures on $\Msw_{[\gamma_{v}]}(\ex B\tc v)$, and so that  the pullback of this Kuranishi structure to $\Msw_{[\gamma]}(\ex B)$ is an embedded Kuranishi structure.

\end{prop}

\pf  
Even though the Kuranishi charts on $\Msw_{[\gamma_{v}]}(\ex B\tc v)$ are compatible with each other, the product Kuranishi charts on $\prod_{v}\Msw_{[\gamma_{v}]}(\ex B\tc v)$  and the pulled-back Kuranishi charts on $\Msw_{[\gamma]}(\ex B)$ may not be compatible (because the product of totally ordered sets usually only has a partial order.) To overcome this problem, we  choose  where to use each chart carefully. The construction is  analogous to the construction of compatible Kuranishi charts in the proof of Theorem 7.3 from \cite{evc}.

In particular, Lemma 7.1 of \cite{evc} and Lemma \ref{sg pullback} above, imply that each holomorphic curve $f$ in $\prod_{v}\Msw_{[\gamma_{v}]}$ has a neighborhood $\mathcal O$ with a $\C\infty 1$ function $\rho\co  \mathcal O\longrightarrow [0,1]$ so that 

\begin{itemize}
\item all holomorphic curves in the closure of $\{\rho>0\}\subset\prod_{v}\Msw_{[\gamma_{v}]}$ are contained in $\mathcal O$;
\item $\rho(f)=1$
\item there exists some collection $(\mathcal U_{v},V_{v},\hat f_{v}/G_{v})$ of the given Kuranishi charts on $\Msw_{[\gamma_{v}]}$ so that $\mathcal O \subset \prod_{v}\mathcal U_{v}$; 
\item $\cut^{*}\oplus_{v}V_{v}$ is simply-generated on $\cut^{-1}(\mathcal O)$.
\end{itemize}

As the set of holomorphic curves in any connected component of $\prod_{v}\Msw_{[\gamma_{v}]}(\ex B\tc v)$ is compact\footnote{To ensure compactness, we have assumed that there is a taming form taming the almost complex structure on $\ex B$, and that $\totb{\ex B}$ admits an immersion into some $\mathbb R^{N}$. These conditions, and the required compactness, then automatically hold for $\ex B\tc v$.}, there exists some finite collection  $\{(\mathcal O_{1},\rho_{1}),\dotsc, (\mathcal O_{N},\rho_{N})\}$ with the substacks $\{\rho_{i}>\frac 12\}$ covering the holomorphic curves in this component of $\prod_{v}\Msw_{[\gamma_{v}]}(\ex B\tc v)$. Use the notation $V_{i}$ for the obstruction bundle, $\oplus_{v}V_{v}$ above,  defined on $\mathcal O_{i}$, and for any $S\subset \{1,\dotsc, N\}$, let $V_{S}:=\oplus_{i\in S} V_{i}$. 

We shall use $V_{S}$ on a restricted domain,   $\mathcal O_{S}$, defined to  be the interior of the following set.
\[\left\{\lrb{\min(0.4,\min_{j\in S}\rho_{j})-\max(0.1,\max_{j'\notin S}\rho_{j'})}>\frac {0.1}N\right\}\subset \bigcap_{j\in S}\mathcal O_{j}\]
In the above, we extend $\rho_{j'}$ to be $0$ wherever it is not already defined. The  proof of Claim 7.4 from \cite{evc} applies without modification to show that that these $\mathcal O_{S}$  form an open cover of the holomorphic curves in our connected component of $\prod_{v}\Msw_{[\gamma_{v}]}(\ex B\tc v)$. Make a similar construction for all other connected components. 

By construction,  $\mathcal O_{S}\cap \mathcal O_{S'}\neq\emptyset$  only if $V_{S}$ is a subsheaf of $V_{S'}$ or visa versa, so we may construct compatible Kuranishi charts on $\prod_{v}\Msw_{[\gamma_{v}]}$ using subcharts of the product chart that uses $V_{S}$ only within $\mathcal O_{S}$.

Let us verify that each of our Kuranishi charts, $(\mathcal O,V,\hat f/G)$ on $\prod_{v}\Msw_{[\gamma_{v}]}(\ex B\tc v)$, pulls back to a Kuranishi chart $(\cut^{-1}\mathcal O,\cut^{*}V,\cut^{*}\hat f/G)$ on $\Msw_{\gamma}(\ex B)$. Our assumption that our original charts were $ev^{[\gamma_{v}]}$--submersive implies that our chart is $(\prod_{v}ev^{[\gamma_{v}]})$--submersive, and in particular, $\prod_{v}ev^{[\gamma_{v}]}$ restricted to $\hat f$ is a submersion. Proposition \ref{cut*} then tells us that $\cut^{*}\hat f$ is a well-defined $\C\infty1$ family of curves, and that $\cut^{*}\hat f/G$ represents $\dbar^{-1}(\cut^{*}V)\subset \cut^{-1}(\mathcal O)$.  We also need that $D\dbar$ is strongly transverse\footnote{See definitions 2.26 and 2.29 of \cite{evc}.} to $\cut^{*}V$ at any holomorphic curve $f$ in $\hat f$. By our submersive assumption, $D\dbar$ at $\cut f$ is strongly transverse to $V$, even when restricted to the kernel of the derivative of $\prod_{v}ev^{[\gamma_{v}]}$. Proposition \ref{cut*} then implies that $D\dbar$ is strongly transverse to $\cut^{*}V$. We already know that $\cut^{*}V$ is simply-generated, and $\cut^{*}V$ is complex because $V$ is, so $(\cut^{-1}\mathcal O,\cut^{*}V,\cut^{*}\hat f/G)$ satisfies all the requirements to be a Kuranishi chart from Definition 2.27 of \cite{evc}.

The pullback of any compatible collection of Kuranishi charts is compatible, so our embedded Kuranishi structure on $\prod_{v}\Msw_{[\gamma_{v}]}(\ex B\tc v)$ pulls back to an embedded Kuranishi structure on $\Msw_{[\gamma]}(\ex B)$.

\stop

Proposition \ref{cut*} almost describes $\cut^{*}\hat f$ as a fiber-product. The following lemma proves that after applying tropical completion as in section 7 of \cite{vfc}, we get an honest fiber-product.

\begin{lemma}\label{tcfp}Given any Kuranishi chart $(\mathcal U,V,\hat f/G)$ on $\prod_{v}\Msw_{[\gamma_{v}]}(\ex B\tc v)$ that pulls back to a Kuranishi chart $(\cut^{-1}U,\cut^{*}V,\cut^{*}\hat f/G)$,  the following is a fiber-product diagram.
\[\begin{tikzcd}\ex F(\cut^{*}\hat f)\tc\gamma\dar\rar{\ev^{[\gamma]}\tc\gamma}&(\End^{[\gamma]}\ex B)\tc {\totb{\ev^{[\gamma]}}\gamma}\dar{\Delta}
\\ \ex F(\hat f)\rar{\prod_{v}\ev^{[\gamma_{v}]}} & \prod_{v}\End^{[\gamma_{v}]}(\ex B\tc v) \end{tikzcd}\]
\end{lemma}

As first glance, the above diagram requires tropical completion of the bottom row at $\coprod_{v}\gamma_{v}$ to make sense, however as we shall see in the proof, such tropical completion does nothing.

\

\pf 

Let us describe the fiber-product $\ex F(\hat f)\times_{\prod_{v}\End^{[\gamma_{v}]}(\ex B\tc v)}\End^{[\gamma]}(\ex B)$. The image of $\prod_{v}\ev^{[\gamma_{v}]}$ is encoded by the following diagram 
\[\begin{tikzcd} \ex F(\hat f^{+})\ar[bend left]{rr}{\prod_{v}\tilde {\ev}^{[\gamma]}_{0}}\rar[hook]\ar{dr}  &\tilde E_{0}\dar \rar & \prod_{v}\gend^{[\gamma_{v}]}_{0}(\ex B\tc v) 
\\ & \ex F(\hat f)\end{tikzcd}\]
where $\tilde E_{0}$ is the unique extension of the bundle $\ex F(\hat f^{+})\longrightarrow\ex F(\hat f)$ so that the following is a fiber-product diagram.
\begin{equation}\label{Epb}\begin{tikzcd}\tilde E_{1}:=\tilde E_{0}\times_{\ex F(\hat f)} \tilde E_{0}\dar[shift left]\dar[shift right]\rar&\prod_{v}\gend^{[\gamma_{v}]}_{1}(\ex B\tc v)\dar[shift left]\dar[shift right]
\\ \tilde E_{0} \rar & \prod_{v}\gend^{[\gamma_{v}]}_{0}(\ex B\tc v) 
 \end{tikzcd}\end{equation}
 Importantly, $\ex F(\hat f^{+})$ is a subset of $\tilde E_{0}$ determined by an open condition on the tropical part of $\tilde E_{0}$. 
 
 \begin{claim}\label{qsc} The fiber-product 
 \[\ex F(\hat f)\times_{\prod_{v}\End^{[\gamma_{v}]}(\ex B\tc v)}\End^{[\gamma]}(\ex B)\]
 is the quotient stack of the groupoid
 \[\tilde E\times_{\prod_{v}\gend^{[\gamma_{v}]}(\ex B\tc v)}\gend^{[\gamma]}(\ex B)\ .\]
 \end{claim}
 
 To prove Claim \ref{qsc}, consider a map of an exploded manifold $\ex X$ into the quotient stack of the above fiber-product of groupoids. Such a map is a bundle $\tilde{\ex X}_{0}\longrightarrow \ex X$ and a pullback diagram.
  \[\begin{tikzcd}\tilde {\ex X}_{1}:=\tilde{\ex X}_{0}\times_{\ex X}\tilde{\ex X}_{0}\dar[shift left]\dar[shift right] \rar&\tilde E_{1}\times_{\prod_{v}\gend_{1}^{[\gamma_{v}]}(\ex B\tc v)}\gend_{1}^{[\gamma]}(\ex B)\dar[shift left]\dar[shift right]
\\ \tilde{\ex X}_{0}\rar& \tilde E_{0}\times_{\prod_{v}\gend_{0}^{[\gamma_{v}]}(\ex B\tc v)}\gend_{0}^{[\gamma]}(\ex B)  \end{tikzcd}\]
 Because of the pullback diagram (\ref{Epb}), the above is a pullback diagram if and only if the induced diagram
 \begin{equation}\begin{tikzcd}\label{XEd}\tilde{\ex X}_{1}\dar[shift left]\dar[shift right]\rar&\gend^{[\gamma]}_{1}(\ex B)\dar[shift left]\dar[shift right]
 \\ \tilde{\ex X}_{0}\rar &\gend^{[\gamma]}_{0}(\ex B) \end{tikzcd}\end{equation}
is a pullback diagram ---  such a diagram is a map of $\ex X$ to $\End^{[\gamma]}(\ex B)$.

 A map of $\ex X$  to our stack fiber-product
 is a map $h\co  \ex X\longrightarrow \ex F(\hat f)$, a diagram in the form (\ref{XEd}), and an isomorphism between their images in $\prod_{v}\End^{[\gamma_{v}]}(\ex B\tc v)$. Such an isomorphism  amounts to a lift of $h$ to a groupoid map $\tilde h\co  \tilde {\ex X}\longrightarrow \tilde E$ so that the following diagram commutes.
 \[\begin{tikzcd}\tilde {\ex X}\dar{\tilde h} \rar&\gend^{[\gamma]}(\ex B)\dar
 \\ \tilde E\rar & \prod_{v}\gend^{[\gamma_{v}]}(\ex B\tc v)\end{tikzcd}\]
 This data of our map of $\ex X$ into the stack fiber-product is equivalent to a map  of $\ex X$ into the quotient stack of our fiber-product of groupoids, so Claim \ref{qsc} is true. 
 
 We need to relate $\ex F(\cut^{*}\hat f)$ to this fiber-product. Proposition \ref{cut*} gives ${\tilde {\ex F}}(\cut^{*}\hat f)$ as a fiber-product, however this fiber-product involves $\hat f'\subset \hat f$, the pullback of $\hat f$ under the inclusion $\prod_{v}\Msw_{[\gamma_{v}]}(\ex B)\longrightarrow \prod_{v}\Msw_{[\gamma_{v}]}(\ex B\tc v)$. This family $\hat f'$ is the subfamily of $\hat f$ consisting of curves with tropical parts having all vertices contained in $\totb{\bar{\ex B}_{v}}\subset\ex B\tc v$. In particular, $\ex F(\hat f')\subset \ex F(\hat f)$ is a subset determined by restricting  to the inverse image of an open\footnote{Note that open subsets of $\ex B$  have closed image in   $\totb{\ex B}$.} subset of $\totb{\ex F(\hat f)}$. Therefore, Proposition \ref{cut*} gives that $\ex F(\hat f)$ represents the quotient stack of a subgroupoid of $\tilde E\times_{\prod_{v}\gend^{[\gamma_{v}]}(\ex B\tc v)}\gend^{[\gamma]}(\ex B)$ determined by restricting to the inverse image of an open subset of its tropical part. It follows that $\ex F(\cut^{*}\hat f)$ represents a subset of the fiber-product \[\ex F(\hat f)\times_{\prod_{v}\End^{[\gamma_{v}]}(\ex B\tc v)}\End^{[\gamma]}(\ex B)\] determined by restricting to the inverse image of an open subset of its tropical part. This subset includes all the points corresponding to curves with tropical part actually equal to $\coprod_{v}\gamma_{v}$, and points in $\End^{[\gamma]}(\ex B)$ that are the image of curves with tropical part $\gamma$, so we may describe the subset, $\ex F(\cut^{*}\hat f)\rvert_{\gamma}$, corresponding to curves with tropical part equal to $\gamma$ as an honest fiber-product.  Taking tropical completions\footnote{Section 7 of \cite{vfc} only describes tropical completion of orbifolds, not general stacks, so we must specify what is meant by tropical completion of $\End^{[\gamma]}\ex B$. This stack is a global quotient of $\rend^{[\gamma]}\ex B$, so define its tropical completion as the corresponding quotient of  the tropical completion of $\rend^{[\gamma]}\ex B$.} at the relevant points corresponding to $\gamma$ gives a fiber-product diagram.
 \[\begin{tikzcd}[column sep=large]\ex F(\cut^{*}\hat f)\tc \gamma \rar{\ev^{[\gamma]}\tc{\gamma}}\dar & (\End^{[\gamma]}\ex B)\tc {\totb{\ev^{[\gamma]}}\gamma}\dar{\Delta\tc{\totb{\ev^{[\gamma]}}\gamma}}
 \\\ex F(\hat f)\tc{\coprod_{v}\gamma_{v}}\rar[swap]{\prod_{v}\ev^{[\gamma_{v}]}\tc{\gamma_{v}}} &\prod_{v}(\End^{[\gamma_{v}]}(\ex B\tc v))\tc{\totb{\ev^{[\gamma_{v}]}}\gamma_{v}}\end{tikzcd}\]
The tropical completions on the bottom row do nothing, because the spaces involved are already complete, and have tropical parts which are always an infinite cone around the tropical completion point. This is because $\totb{\ex B\tc v}$ is an infinite cone around $v$. The rescaling action around $v$ also acts on $\totb{\End^{[\gamma_{v}]}(\ex B\tc v)}$, preserving the image of curves with tropical part $\gamma_{v}$ (and only this point), so tropical completion at this point does nothing. Similarly, because $\hat f$ has universal tropical structure, the rescaling action also acts on $\totb{\ex F(\hat f)}$. This scaling action preserves only the point corresponding to curves with tropical part $\gamma_{v}$, so again tropical completion at this point does nothing. 

Removing the unnecessary tropical completions from the bottom row gives the required fiber-product diagram.

\stop

We are now ready to write our first gluing theorem. Consider the   map $ev\co  \Msw_{\cdot}(\ex B)\longrightarrow \coprod_{n}(\rend\ex B)^{n}$, and define 
\[\eta:= ev_{!}(\hbar^{2g-2+n}q^{E})\]
 where $g$, $n$ and $E$ are the locally constant functions recording the genus, number of ends,  and $\omega$--energy of curves, and $\hbar$ and $q$ are dummy variables.
In general, $ev^{[\gamma]}$ is not complete, even restricted to curves with bounded genus and energy, and $ev^{[\gamma]}_{!}$ may not be defined. Use tropical completion as in section 7 of \cite{vfc} to define
 \[\eta^{[\gamma]}:=(ev^{[\gamma]}\tc\gamma)_{!}(\hbar^{2g-2+n}q^{E})\]
 and similarly use tropical completion at $\gamma$ to define
 \[\eta\tc\gamma:=(ev\tc\gamma)_{!}(\hbar^{2g-2+n}q^{E})\ .\]
 Similarly, given any complex vector-bundle $W$ over $\Msw_{\bullet}(\ex B)$, define
 \[\eta(W):=ev_{!}(\hbar^{2g-2+n}q^{E}c(W))\]
 \[\eta^{[\gamma]}(W):=(ev^{[\gamma]}\tc\gamma)_{!}(\hbar^{2g-2+n}q^{E}c(W\tc\gamma))\]
 \[\eta(W)\tc\gamma:=(ev\tc\gamma)_{!}(\hbar^{2g-2+n}q^{E}c(W\tc\gamma))\]
 where $c(W)$ is the top Chern class of $W$.

\begin{remark} In the case of $\gamma_{v}$,  tropical completion at $\gamma_{v}$ affects nothing, and $\eta^{[\gamma_{v}]}(W)$ is the the restriction of $\eta(W)$ to $\rend^{[\gamma_{v}]}(\ex B\tc v)\subset \coprod_{n}(\rend \ex B\tc v)^{n}$.
 \end{remark}
\begin{remark}Given any point $p\in \coprod_{n}(\totb{\rend\ex B})^{n}$, say that $\gamma\in \totb{ev}^{-1}p$ if  curves with tropical part $\gamma$ are send by $ev$ to points with tropical part $p$. Lemma 7.7 of \cite{vfc} gives that
\[\eta\tc p=\sum_{\gamma\in \totb{ev}^{-1}p}\eta\tc \gamma\]
and more generally, 
\[\eta(W)\tc p=\sum_{\gamma\in \totb{ev}^{-1}p}\eta(W)\tc \gamma\ .\]
\end{remark}

  Note that the target of $ev^{[\gamma]}\tc \gamma$ is not $\rend^{[\gamma]}(\ex B)$, but its tropical completion at the image of $\gamma$, so $\eta^{[\gamma]}$ is a (refined\footnote{The minimal cohomology theory of exploded manifolds containing the usual cohomology but with pushforwards compatible with fiber-products is called refined cohomology. See section 9 of \cite{dre}.}) cohomology class on $\rend^{[\gamma]}(\ex B)\tc {\totb{ev^{[\gamma]}}(\gamma)}$.
Indicate the tropical completion of $\Delta$ again by  $\Delta$.
\[\Delta\co  \rend^{[\gamma]}(\ex B)\tc {\totb{ev^{[\gamma]}}(\gamma)}\longrightarrow \prod_{v}\rend^{[\gamma_{v}]}(\ex B\tc v)\]

 \begin{thm}\label{sg}
\[ \Delta^{*}\prod_{v}\eta^{[\gamma_{v}]}=k_{\gamma}\eta^{[\gamma]}\]
where $k_{\gamma}=\prod_{e\in\ie\gamma}m_{e}$, the product of the multiplicities of the internal edges of $\gamma$.
Similarly, given  complex vector bundles $W_{v}$ on $\Msw_{[\gamma_{v}]}$, 
\[ \Delta^{*}\prod_{v}\eta^{[\gamma_{v}]}(W_{v})=k_{\gamma}\eta^{[\gamma]}(\cut^{*}\oplus_{v}W_{v})\tc\gamma\ .\]
 \end{thm}
 
 \pf
 
 Using \cite{evc}, construct $ev^{[\gamma_{v}]}$--submersive embedded Kuranishi structures on $\Msw_{[\gamma_{v}]}(\ex B\tc v)$, then use Proposition \ref{Kuranishi fiber-product} to construct a corresponding embedded Kuranishi structure on $\prod_{v}\Msw_{[\gamma_{v}]}(\ex B\tc v)$. The resulting Kuranishi category $\mathcal K$ within $\prod_{v}\Msw_{[\gamma_{v}]}(\ex B\tc v)$ is a weak product of the Kuranishi categories $\mathcal K_{v}$ from $\Msw_{[\gamma_{v}]}(\ex B\tc v)$, so Theorem 6.2 from \cite{vfc} gives the expected product relation  when pushing forward using $\mathcal K$ and $\prod_{v}ev^{[\gamma_{v}]}$, or $\mathcal K_{v}$ and $ev^{[\gamma_{v}]}$. 
 \begin{equation}\label{g3}(\prod_{v}ev^{[\gamma_{v}]})_{!}(\prod_{v}\hbar^{2g_{v}-2+n_{v}}q^{E_{v}}c(W_{v}))=\prod_{v}\eta^{[\gamma_{v}]}(W_{v})\end{equation}
 Now consider the Kuranishi category, $\cut^{*}\mathcal K$, defined using the pullback of our embedded Kuranishi structure, and apply tropical completion\footnote{See section 7 of \cite{vfc}, especially Lemma 7.7.} to $ev^{[\gamma]}$ considered as a map from  $\cut^{*}\mathcal K$.  Lemma \ref{tcfp} implies that the following is a pullback diagram of Kuranishi categories.
 \[\begin{tikzcd}\cut^{*}\mathcal K\tc\gamma \dar{\cut\tc \gamma}\rar{\ev^{[\gamma]}\tc\gamma}&(\End^{[\gamma]}(\ex B))\tc{\totb{\ev^{[\gamma]}} \gamma}\dar{\Delta}
 \\ \mathcal K\rar{\prod_{v}\ev^{[\gamma_{v}]}}& \prod_{v}\End^{[\gamma_{v}]}(\ex B\tc v)\end{tikzcd}\]
 So long as the stacks above on the right are orbifolds, we can apply Theorem 5.22 of \cite{vfc} to obtain the relationship between $(\ev^{[\gamma]}\tc\gamma)_{!}$ and $(\prod_{v}\ev^{[\gamma_{v}]})_{!}$. We have that $\End^{[\gamma]}(\ex B)=\rend^{[\gamma]}(\ex B)/\prod_{e}G_{e}$, where the group $G_{e}$ is $\mathbb Z_{m_{e}}$ for each edge of multiplicity $m_{e}\neq 0$, and is $\ex T$ for each edge of multiplicity $0$.
 We can forget the $G_{e}$--action for each end $e$ of $\gamma$, leaving us with the following pullback diagram.
 \[\begin{tikzcd}\cut^{*}\mathcal K\tc\gamma \dar{\cut\tc \gamma}\rar&(\rend^{[\gamma]}(\ex B))\tc{\totb{ev^{[\gamma]}} \gamma}/\prod_{e\in\ie\gamma }G_{e}\dar
 \\ \mathcal K\rar& \prod_{v}\rend^{[\gamma_{v}]}(\ex B\tc v)/\prod_{e\in\ie\gamma}G_{e}^{2}\end{tikzcd}\]
In particular, $\cut^{*}\mathcal K\tc \gamma$ is a $(\prod_{e\in\ie\gamma}G_{e})$--bundle over the corresponding fiber-product forgetting the $G_{e}$--actions. If any internal edge has multiplicity $0$, it follows that $\eta^{[\gamma]}$ must be $0$.  Otherwise, the righthand side of the above consists of orbifolds, so applying Theorem 5.22 of \cite{vfc}, then pushing forward the result via the map forgetting the $G_{e}$--actions, gives 
\begin{equation}\label{g1}(ev^{[\gamma]}\tc\gamma)_{!}(\cut^{*}\theta)=\lrb{\prod_{e\in\ie\gamma}m_{e}}\Delta^{*}(\prod_{v}ev^{[\gamma_{v}]})_{!}(\theta)\end{equation}
 Observing that $\cut^{*}(\prod_{v}\hbar^{2g_{v}-2+n_{v}}q^{E_{v}}c(W_{v}))=\hbar^{2g-2+n}q^{E}c(\cut^{*}\oplus_{v}W_{v})$, then combining equations (\ref{g1}), and (\ref{g3}) gives our desired result.
 \[ \Delta^{*}\prod_{v}\eta^{[\gamma_{v}]}(W_{v})=k_{\gamma}\eta^{[\gamma]}(\cut^{*}\oplus_{v}W_{v})\]

 \stop
 
 Consider the map
\[i^{[\gamma]}\co  \rend^{[\gamma]}(\ex B)\tc{\totb{ev^{[\gamma]}}\gamma}\longrightarrow \lrb{\coprod_{n}(\rend \ex B)^{n}}\tc{\totb{ev}\gamma}\] forgetting all internal edges.

 \begin{lemma}\label{eta pf} The following relationship holds between $\eta\tc \gamma$ and $\eta^{[\gamma]}$.
 \[\eta\tc \gamma=\frac 1{\abs{\Aut\gamma}}i^{[\gamma]}_{!}\eta^{[\gamma]}\]
More generally, let $W$ be a complex vectorbundle on $\Msw_{\cdot}$, and let $\pi^{*}W$ indicate its pullback under the map $\pi\co  \Msw_{[\gamma]}\longrightarrow \Msw_{\cdot}$. The following equation holds.
\[\eta(W)\tc \gamma=\frac 1{\abs{\Aut\gamma}}i^{[\gamma]}_{!}\eta^{[\gamma]}(\pi^{*}W)\]
 \end{lemma}
 
 \pf
 
 Choose an embedded Kuranishi structure on $\Msw_{\cdot}(\ex B)$ for defining $\eta$, and pull back this embedded Kuranishi structure using $\pi$ to define an embedded Kuranishi structure on $\Msw_{[\gamma]}(\ex B)$. Let $\mathcal K$ and $\pi^{*}\mathcal K$ be the associated Kuranishi categories. Take the tropical completion of these Kuranishi categories at $\gamma$, and consider the following commutative diagram.
 
 \[\begin{tikzcd}\mathcal K\tc \gamma\rar{ev\tc \gamma}&\lrb{\coprod_{n} (\rend \ex B)^{n}}\tc {\totb{ev}\gamma}
 \\ \pi^{*}\mathcal K\tc \gamma\uar{\pi}\rar{ev^{[\gamma]}}&\lrb{\rend^{[\gamma]}\ex B}\tc{\totb{ev^{[\gamma]}}\gamma}\uar{i^{[\gamma]}}
 \end{tikzcd}\]
 
 When $\mathcal K$ uses the family $\hat f$ with automorphism group $G$,  $\pi^{*}\mathcal K$ uses  $\pi^{*}\hat f$ with the action of $G$ induced using the universal property of $\pi^{*}\hat f$. For a given curve $f$ with tropical part $\gamma$, there are $\abs{\Aut\gamma}$ ways of $\gamma$--decorating $f$ so that the $\gamma$--decoration is an isomorphism. It follows that $\ex F(\pi^{*}\hat f)\tc \gamma\longrightarrow \ex F(\hat f)\tc\gamma$ is an $\abs{\Aut\gamma}$--fold cover, so $\pi\co  \pi^{*}\mathcal K\tc \gamma\longrightarrow \mathcal K\tc \gamma$ is an $\abs{\Aut\gamma}$--fold cover. The required formula for pushforwards follows.
 
 \stop
 
 \section{Enhanced evaluation map and gluing formula}\label{ee map}

In this section, we enhance our evaluation map to generalize the stabilization map $ev^{0}\co  \Msw_{g,n}(\ex B)\longrightarrow \Msw_{g,n}(pt)$, where $\Msw_{g,n}(pt)$ is the moduli stack of stable exploded curves with genus $g$ and $n$ marked points.  This stabilization map is constructed in section 4.1 of \cite{evc}, where it is also shown that $\Msw_{g,n}(pt)$ is an orbifold and the explosion of the corresponding Deligne-Mumford space relative to its boundary  divisors. 

For curves with genus $g$ and $n$ punctures where $2g-2+n>0$, the target, $\mathcal X_{g,n}(\ex B)$, of our enhanced evaluation map can be thought of as a fiber-product.
\[\begin{tikzcd}\mathcal X_{g,n}(\ex B)\rar\dar &\Msw_{g,n}(pt)\dar
\\ (\End \ex B)^{n}\rar &\cdot/\ex T^{n}\end{tikzcd}\]
Construct $\mathcal X_{g,n}(\ex B)$ as the quotient stack of a groupoid $\tilde{\mathcal X}_{g,n}(\ex B)$ with objects  as follows.
\[(\tilde{\mathcal X}_{g,n}(\ex B))_{0}=(\gend_{0}\ex B)^{n}\times \Msw_{g,n^{+}}(pt)\]
 In the case $2g-2+n>0$,  define $\Msw_{g,n^{+}}(pt)$ in analogy to $\Msw_{[\gamma]^{+}}$ as the stack of stable curves in $\Msw_{g,n}$ along with an extra choice of point in each of the $n$ ends. When $(g,n)$ is $(0,0)$, $(0,1)$, or $(1,0)$, define $\Msw_{g,n^{+}}(pt)$ to be a point, and in the remaining case, define $\Msw_{0,2^{+}}(pt)$ to be the stack of curves isomorphic to $\ex T$ with an extra choice of $2$ points with distinct tropical part. This $\Msw_{0,2^{+}}(pt)$ is canonically isomorphic to $\et 1{(0,\infty)}$, and to the stratum of $\Msw_{0,4}(pt)$ pairing the 1st and 3rd, and 2nd and 4th ends. 
 
 In each case, there is a (partially defined) action of $\ex T^{n}$ on $\Msw_{g,n^{+}}(pt)$, moving the $n$ extra points. This action corresponds to the action on $\Msw_{[\gamma]^{+}}$ in the stable case. For a curve in $\Msw_{0,2^{+}}$, we may fix an isomorphism with $\ex T$ so that the image of the first extra point in $\totb{\ex T}$ is before the second extra point. Then $(c_{1},c_{2})\in \ex T^{2}$ acts by multiplying the first point by $c_{1}^{-1}$ and the second point by $c_{2}$, and is defined so long as $\totb{c_{1}/c_{2}}$ is less than the distance between the image of our points in $\totb{\ex T}=\mathbb R$. Recalling that there is also a (partially defined) action of $\ex T$ on $\gend_{0}\ex B$, let $\tilde{\mathcal X}_{g,n}(\ex B)$  be the groupoid defined by the (partially defined) action of $\ex T^{n}$ on $(\gend_{0}\ex B)^{n}\times \Msw_{g,n^{+}}(pt)$, so the morphisms in our groupoid $\tilde{\mathcal X}_{g,n}(\ex B)$ are parametrized by the subset
 \[(\tilde{\mathcal X}_{g,n}(\ex B))_{1}\subset (\gend_{1}\ex B)^{n}\times \Msw_{g,n^{+}}(pt)\subset \ex T^{n}\times (\gend_{0}\ex B)^{n}\times \Msw_{g,n^{+}}(pt) \]
where the action of $c\in \ex T^{n}$ on $p\in (\gend_{0}\ex B)^{n}\times \Msw_{g,n^{+}}(pt) $ is defined. The two maps $(\tilde{\mathcal X}_{g,n}(\ex B))_{1}\rightrightarrows (\tilde{\mathcal X}_{g,n}(\ex B))_{0}$ are given by $(c,p)\to p$ and $(c,p)\to c*p$. This defines the groupoid $\tilde {\mathcal X}_{g,n} (\ex B)$. It is also convenient to take the union of these groupoids for all $n$ and $g$. 

 \[ \tilde {\mathcal X}:= \coprod_{g,n}\tilde{\mathcal X}_{g,n} \]  

Our evaluation map $\tilde{\ev}$ extends to the enhanced evaluation map
\[\tilde {EV}\co  \tilde{\Msw_{\cdot}}(\ex B)\longrightarrow \tilde {\mathcal X}(\ex B)\]
with
\[\tilde {EV}_{0}\co  \Msw_{\cdot^{+}}(\ex B)\longrightarrow \tilde{\mathcal X}_{0}(\ex B)\]
given by $\tilde{\ev}_{0}$ on the first factor, and the stabilization map $\Msw_{\cdot^{+}}(\ex B)\longrightarrow \coprod_{g,n}\Msw_{g,n^{+}}$ on the second factor (treating $\Msw_{g,n^{+}}$ as a substack of $\Msw_{g,2n}$). This $\tilde{EV}_{0}$ is  equivariant with respect to the various (partially defined) $\ex T^{n}$--actions, and lifts uniquely to the map of groupoids $\tilde{EV}$ so that the following diagram commutes.
\[\begin{tikzcd}\tilde {\Msw_{\cdot}}(\ex B)\rar{\tilde{EV}}\ar{dr}{\tilde{\ev}}&\tilde{\mathcal X}(\ex B)\dar
\\ & \coprod_{n}(\gend\ex B)^{n}\end{tikzcd}\] 
Letting $\mathcal X$ be the quotient stack of $\tilde{\mathcal X}$, we get the following commutative diagram of evaluation maps:
\[\begin{tikzcd}\Msw_{\cdot}(\ex B)\rar{ev^{0}}\ar{dr}{EV}\dar{\ev}& \coprod_{g,n}\Msw_{g,n}(pt)
\\ \coprod_{n}(\End\ex B)^{n}&\lar\uar \mathcal X(\ex B)\end{tikzcd}\]

In the case that $\totb{\ex B}$ is bounded and $2g-2+n>0$, $\mathcal X_{g,n}=\ex B^{n}\times \Msw_{g,n}(pt)$, and $EV$ is (the exploded version of) a familiar evaluation map used in Gromov--Witten theory. 

As with $\ev$, we can enhance $\ev^{[\gamma]}$ to obtain $EV^{[\gamma]}$. For a vertex $v$ of $\gamma$, use the notation $\Msw_{[\gamma_{v}]^{+}}(pt):=\coprod_{g}\Msw_{g,n^{+}}(pt)$, where we identify the $n$ edges of $\gamma_{v}$ with the $n$ labels from $\Msw_{g,n}(pt)$.  Define
\[\tilde{\mathcal X}^{[\gamma]}_{0}(\ex B):=\gend^{[\gamma]}_{0}\ex B\times \prod_{v}\Msw_{[\gamma_{v}]^{+}}(pt)\ .\]
Then define
\[\tilde {EV}^{[\gamma]}_{0}\co  \Msw_{[\gamma]^{+}}(\ex B)\longrightarrow \tilde {\mathcal X}^{[\gamma]}_{0}(\ex B)\]
as $\tilde \ev^{[\gamma]}_{0}$ on the first factor, and the cutting map followed by the stabilization map on the second factor. $\tilde{EV}^{[\gamma]}_{0}$ is equivariant with respect to the (partially defined) $\ex T$ action corresponding to each edge of $\gamma$, (acting diagonally by $(c,c^{-1})$ on $\prod_{v}\Msw_{[\gamma_{v}]^{+}}(pt)$ in the case of an internal edge of $\gamma$)  so we can promote $\tilde {\mathcal X}^{[\gamma]}_{0}$ to a groupoid and $\tilde {EV}^{[\gamma]}_{0}$ to a groupoid map so that $\tilde {\mathcal X}_{1}^{[\gamma]}(\ex B)$ is a subset of $\gend^{[\gamma]}_{1}(\ex B)\times \prod_{v}\Msw_{[\gamma_{v}]^{+}}(pt)$, and so that the following commutative diagram of groupoid maps exists. 
\[\begin{tikzcd}[column sep=large]&\gend^{[\gamma]}(\ex B)
\\\tilde{\Msw_{[\gamma]}}(\ex B)\rar{\tilde {EV}^{[\gamma]}}\ar{ur}{\tilde{\ev}^{[\gamma]}}\dar{cut}&\tilde{\mathcal X}^{[\gamma]}(\ex B)\uar\dar{\tilde \Delta}
\\ \prod_{v}\tilde{\Msw_{[\gamma_{v}]}}(\ex B\tc v)\rar{\prod_{v}\tilde{EV}^{[\gamma_{v}]}} &\prod_{v}\tilde{\mathcal X}^{[\gamma_{v}]}(\ex B\tc v)  \end{tikzcd}\]
Passing to quotient stacks, we get the following commutative diagram
\[\begin{tikzcd}[column sep=large]&\End^{[\gamma]}(\ex B)
\\{\Msw_{[\gamma]}}(\ex B)\rar{  {EV}^{[\gamma]}}\ar{ur}{ {\ev}^{[\gamma]}}\dar{cut}& {\mathcal X}^{[\gamma]}(\ex B)\uar\dar{  \Delta}
\\ \prod_{v} \Msw_{[\gamma_{v}]}(\ex B\tc v)\rar{\prod_{v} {EV}^{[\gamma_{v}]}} &\prod_{v} {\mathcal X}^{[\gamma_{v}]}(\ex B\tc v)  \end{tikzcd}\]

For notational convenience, we shall use the same notation for $\Delta$ and its tropical completion.
\[\Delta\co  (\mathcal X^{[\gamma]}\ex B)\tc{\totb{EV^{[\gamma]}}\gamma}\longrightarrow \prod_{v}\mathcal X^{[\gamma_{v}]}(\ex B\tc v)\]

\begin{lemma}\label{xfp}The following is a fiber-product diagram.
\[\begin{tikzcd}(\mathcal X^{[\gamma]}\ex B)\tc{\totb{EV^{[\gamma]}}\gamma}\rar\dar{\Delta}&\End^{[\gamma]}(\ex B)\tc{\totb{\ev^{[\gamma]}}\gamma}\dar{\Delta}
\\ \prod_{v}\mathcal X^{[\gamma_{v}]}(\ex B\tc v) \rar & \prod_{v}\End^{[\gamma_{v}]}(\ex B\tc v)\end{tikzcd}\]
\end{lemma}

\pf

Applying tropical completion to the above stacks at the image of $\gamma$ is equivalent to applying tropical completion to the corresponding groupoids at the image of curves with tropical part $\gamma$ and extra chosen points some fixed location on $\gamma$. Although tropical completion at the image of $\gamma_{v}$ does nothing to the stacks $\mathcal X^{[\gamma_{v}]}$ and $\End^{[\gamma_{v}]}$, applying tropical completion to the corresponding groupoids has the effect of replacing  our partially defined $\ex T$--actions with honest $\ex T$--actions. After applying tropical completion at appropriate points, we get the following commutative diagram of groupoids
\[\begin{tikzcd}\tilde{\mathcal X}^{[\gamma]}\tc a\rar\dar &\gend^{[\gamma]}(\ex B)\tc b\dar
\\ \prod_{v}\tilde{\mathcal X}^{[\gamma_{v}]}\tc c\rar & \prod_{v}\gend^{[\gamma_{v}]}(\ex B\tc v)\tc d
\end{tikzcd}\]
where all groupoid actions are honest actions of some $\ex T^{n}$. This is a fiber product diagram, because the left side is the product of the right with the tropical completion of $\prod_{v}\Msw_{[\gamma_{v}]^{+}}(pt)$, at both level $0$ and $1$.  Passing to quotient stacks therefore gives the desired fiber-product diagram of stacks.
\[\begin{tikzcd}\mathcal X^{[\gamma]}\tc{\totb{EV^{[\gamma]}}\gamma}\rar\dar&\End^{[\gamma]}(\ex B)\tc{\totb{\ev^{[\gamma]}}\gamma}\dar
\\ \prod_{v}\mathcal X^{[\gamma_{v}]}\rar & \prod_{v}\End^{[\gamma_{v}]}(\ex B\tc v)\end{tikzcd}\]

\stop

In analogy with our definition of $\eta$, define
\[\mu:= EV_{!}(q^{E})\]
\[\mu\tc \gamma:= (EV\tc \gamma)_{!} (q^{E})\]
\[\mu^{[\gamma]}:=(EV^{[\gamma]}\tc\gamma)_{!}(q^{E})\]
where we only need one dummy-variable, $q$, because genus is automatically tracked in $\mathcal X$. These pushforwards are defined as in \cite{vfc} on connected components of $\mathcal X$ that are orbifolds. On the other components, define $\mu$ to vanish.  For $W$ a complex vectorbundle over $\Msw_{\bullet}$, also define
\[\mu:= EV_{!}(q^{E}c(W))\]
\[\mu(W)\tc \gamma:= (EV\tc \gamma)_{!} (q^{E}c(W))\]
\[\mu^{[\gamma]}(W):=(EV^{[\gamma]}\tc\gamma)_{!}(q^{E}c(W))\]
 Note that $\mu^{[\gamma]}$ is a refined cohomology class on $\mathcal X^{[\gamma]}\tc {\totb{EV^{[\gamma]}}\gamma}$, however $\mathcal X^{[\gamma_{v}]}\tc{\totb{EV^{[\gamma_{v}]}}\gamma_{v}}$ coincides with $\mathcal X^{[\gamma_{v}]}$, and $\mu^{[\gamma_{v}]}$ is the restriction of $\mu$ to $\mathcal X^{[\gamma_{{v}}]}\subset\mathcal X(\ex B\tc v)$.

\begin{thm}\label{xg}
\[\mu^{[\gamma]}=\Delta^{*}\prod_{v}\mu^{[\gamma_{v}]}\]
and given complex vector bundles $W_{v}$ on $\Msw_{[\gamma_{v}]}(\ex B\tc v)$, 
\[\mu^{[\gamma]}(\cut^{*}\oplus_{v}W_{v})=\Delta^{*}\prod_{v}\mu^{[\gamma_{v}]}(W_{v})\ .\]
\end{thm}

\pf
As in the proof of Theorem \ref{sg}, choose $ev^{[\gamma_{v}]}$--submersive embedded Kuranishi structures on $\Msw_{[\gamma_{v}]}(\ex B\tc v)$, then use Proposition \ref{Kuranishi fiber-product} to construct a corresponding embedded Kuranishi structure on $\prod_{v}\Msw_{[\gamma_{v}]}(\ex B\tc v)$ that pulls back to an embedded Kuranishi structure using $\cut^{*}$. Let $\mathcal K$ and $\cut^{*}\mathcal K$ be the corresponding Kuranishi categories. Lemmas \ref{xfp} and \ref{tcfp} imply that the following is a fiber-product diagram.
\[\begin{tikzcd}\cut^{*}\mathcal K\tc \gamma\rar{EV^{[\gamma]}\tc \gamma}\dar{\cut}& \mathcal X^{[\gamma]}\tc{\totb{EV} \gamma}\dar{\Delta}\\ 
\mathcal K \rar{\prod_{v}EV^{[\gamma_{v}]}} &\prod_{v}\mathcal X^{[\gamma_{v}]} \end{tikzcd}\]

The only non-orbifold components of $\mathcal X^{[\gamma_{v}]}$ concern curves with genus $0$, and one or two ends, all with zero multiplicity.  The non-orbifold components of $\mathcal X^{[\gamma]}\tc{\totb{EV} \gamma}$ are the inverse image of  non-orbifold components of $\prod_{v}\mathcal  X^{[\gamma_{v}]}$. In the remaining orbifold cases, Theorem 5.22 of \cite{vfc} applies, and our desired formula follows from Theorems 5.22 and 6.2 of \cite{vfc}, as in the proof of Theorem \ref{sg}.

\stop

\

Theorem \ref{stacks fiber-product} allows us to glue together a family of  curves in $\prod_{v}\Mpsw_{[\gamma_{v}^{+}]}(pt)$ to obtain a family of curves in $\Mpsw_{[\gamma]^{+}}(pt)$.  Forgetting the $\gamma$--decoration and extra points on internal edges gives a family in $\Mpsw_{\cdot^{+}}$. Composing this gluing map with the stabilization map $\Mpsw_{\cdot^{+}}(pt)\longrightarrow \Msw_{\cdot^{+}}(pt)$,   and pre-composing with the inclusion $\Msw_{[\gamma_{v}]^{+}}(pt)\longrightarrow \Mpsw_{[\gamma_{v}^{+}]}(pt)$ gives a map,
\[\prod_{v}\Msw_{[\gamma_{v}]^{+}}(pt)\longrightarrow \coprod_{g,n}\Msw_{g,n^{+}}(pt)\]
compatible with the various (partially defined) actions of $\ex T$, and in particular, invariant under the diagonal action corresponding to each internal edge of $\gamma$. The product of this map with  $\tilde i^{[\gamma]}\co  \gend^{[\gamma]}(\ex B)\longrightarrow (\gend\ex B)^{n}$ defines a map of groupoids 
\[\tilde I^{[\gamma]}\co  \tilde {\mathcal X}^{[\gamma]}\longrightarrow \tilde {\mathcal X}\]
compatible with our evaluation maps. Passing to quotient stacks and taking tropical completion at the image of $\gamma$ gives a map
\[I^{[\gamma]}\co  \mathcal X^{[\gamma]}\tc{\totb{EV^{[\gamma]}}\gamma}\longrightarrow \mathcal X\tc{\totb{EV}\gamma}
\]

\begin{lemma}\label{xpf}The following relationship holds between $\mu$ and $\mu^{[\gamma]}$
\[\mu\tc \gamma=\frac 1{\abs{\Aut\gamma}}I^{[\gamma]}_{!}\mu^{[\gamma]} \]
(where $I^{[\gamma]}_{!}\mu^{[\gamma]}$ is defined to be zero on non-orbifold components of $\mathcal X\tc \gamma$).
More generally, let $W$ be a complex vectorbundle on $\Msw_{\cdot}$ and let $\pi^{*}W$ indicate its pullback under the map $\pi\co  \Msw_{[\gamma]}\longrightarrow \Msw_{\cdot}$. Then
\[\mu(W)\tc \gamma=\frac 1{\abs{\Aut\gamma}}I^{[\gamma]}_{!}\mu^{[\gamma]}(\pi^{*}W)\]
\end{lemma}
\pf

Choose embedded Kuranishi structures and use notation as in the proof of Lemma \ref{eta pf}. In this case, we must consider the following commutative diagram.
\[\begin{tikzcd}\mathcal K\tc \gamma\rar{EV\tc \gamma}&\mathcal X\tc{\totb{EV} \gamma}
\\ \pi^{*}\mathcal K\tc \gamma\uar{\pi}\rar{EV^{[\gamma]}\tc \gamma}&\mathcal X^{[\gamma]}\tc{\totb{EV^{[\gamma]}}\gamma}\uar{I^{[\gamma]}}\end{tikzcd}\]
 As in the proof of Lemma \ref{eta pf}, $\pi$ is an $\abs{\Aut \gamma}$--fold cover. If all components of the righthand side of the above diagram were orbifolds, our formula would follow immediately. We need to check that components of $\pi^{*}\mathcal K\tc\gamma$ sent to non-orbifold components of $\mathcal X^{[\gamma]}\tc{\totb{EV^{[\gamma]}} \gamma}$ do not contribute anything to $\mu(W)\tc\gamma$. If a curve with tropical part $\gamma$ is sent to a non-orbifold component of $\mathcal X^{[\gamma]}$, some vertex $v$ of $\gamma$ must have one or two edges, all with $0$ multiplicity, and the stratum of our curve labelled by $v$ must have genus $0$. Within $\pi^{*}\mathcal K$, such a curve must be contained in a family of curves that allows all possible lengths for the corresponding edges. In each case, there is a $\et 1{(0,\infty)}$--worth of choice for these edge lengths that is crunched to a single point under $I^{[\gamma]}\circ EV^{[\gamma]}\tc \gamma$, therefore the corresponding connected component of $\mathcal K\tc \gamma$ does not contribute to $\mu(W)$. 
 
 As all connected components of $\pi^{*}\mathcal K\tc \gamma$ that contribute to $\mu(W)$ are sent to orbifold components of $\mathcal X^{[\gamma]}\tc {\totb{EV^{[\gamma]}}\gamma}$, we may use the formula $(EV\tc \gamma\circ \pi)_{!}=I_{!}^{[\gamma]}\circ (EV^{[\gamma]}\tc \gamma)_{!}$, and the fact that $\pi$ is an $\abs{\Aut\gamma}$--fold cover to prove the desired  relationship.
 \[\mu(W)\tc \gamma=\frac 1{\abs{\Aut\gamma}}I^{[\gamma]}_{!}\mu^{[\gamma]}(\pi^{*}W)\]
 
 \stop
 
 \
 
 For curves in a smooth symplectic manifold $B$, Lemma \ref{xpf} and Theorem \ref{xg} recover Kontsevich and Mannin's splitting and genus-reduction axioms of Gromov--Witten invariants. For example, the genus-reduction axiom may be understood as follows:  The tropical part of $B$ is a single point. The tropical part of an exploded curve corresponding to a curve in $B$ with a single, non-separating node, genus $g$ and $n$ punctures is a tropical curve $\gamma$ with a single vertex, $n$ ends corresponding to punctures and a single interior edge corresponding to the node.  After fixing the labeling of the $n$ ends, there are $2$ automorphisms of $\gamma$, the nontrivial one reversing the interior edge. 

The relevant component of $\mathcal X(B)$ is $\mathcal X_{g,n}(B)= B^{n}\times \mathcal M_{g,n}$, where $\mathcal M_{g,n}$ is the moduli stack of stable exploded curves with genus $g$ and $n$ labeled ends. This may also be thought of as the explosion of the corresponding Delign-Mumford stack, $\bar M_{g,n}$; see section 4.1 of \cite{evc}. Similarly, the relevant component of $\mathcal X^{[\gamma_{v}]}(B\tc v)$ is $B^{n+2}\times\mathcal M_{g-1,n+2}$, and the genus $g$ component of  $\mathcal X^{[\gamma]}$ is a $\et 1{(0,\infty)}$--bundle over $B^{n+1}\times \mathcal M_{g-1,n+2}$. This extra $\et 1{(0,\infty)}$ bundle is created as follows. Take the $(\et 1{(0,\infty)})^{2}$--bundle over $\mathcal M_{g-1,n+2}$ given by choosing an extra point in each of the last two edges, and quotient this bundle by the (partially defined) diagonal action of $\ex T$ multiplying by $(c,c^{-1})$. This $\et 1{(0,\infty)}$--bundle over $\mathcal M_{g-1,n-2}$ is also a $\mathbb Z_{2}$--fold cover of the stratum of $\mathcal M_{g,n}$ corresponding to curves with a single vertex and single internal edge. The $\mathbb Z_{2}$--action swaps the labels of the last two edges in  $\mathcal M_{g-1,n+2}$. In terms of Deligne-Mumford space, the $\mathbb Z_{2}$--quotient of this $\et 1{(0,\infty)}$--bundle corresponds to the $\mathbb C^{*}$--bundle obtained by removing the zero-section of the normal bundle of $\bar M_{g-1,n+2}/\mathbb Z_{2}\subset \bar M_{g,n}$. When we apply tropical completion, we replace this $\et 1{(0,\infty)}$--bundle with a $\ex T$--bundle, $\ex T\rtimes \mathcal M_{g-1,n+2}$. 

Our gluing formula is
\[\mu\tc \gamma=\frac 12I^{[\gamma]}_{!}\Delta^{*}\mu^{[\gamma_{v}]}\]
 stated in terms of the maps 
\[ B^{n}\times B^{2}\times \mathcal M_{g-1,n+2}\xleftarrow{\Delta} B^{n}\times B\times \ex T\rtimes \mathcal M_{g-1,n+2}\xrightarrow{I^{[\gamma]}} B^{n}\times (\ex T\rtimes \mathcal M_{g-1,n+2})/\mathbb Z_{2}\]
where $\Delta$ is the product of the identity on $B^{n}$ with the diagonal $B\longrightarrow B^{2}$ and the bundle map 
$\ex T\rtimes \mathcal M_{g-1,n+2}\longrightarrow \mathcal M_{g-1,n+2}$, and $I^{[\gamma]}$ is the composition of a projection, $\pi$, crushing the extra factor of $\ex B$, and a $\mathbb Z_{2}$--fold covering map,
\[\psi\co  \ex B^{n}\times (\ex T\rtimes \mathcal M_{g-1,n+2})\longrightarrow \ex B^{n}\times (\ex T\rtimes \mathcal M_{g-1,n+2})/\mathbb Z_{2}\ .\]
 Our gluing formula may be rewritten as follows. \begin{equation}\label{kmgf}\psi^{*}\mu\tc \gamma=\pi_{!}\Delta^{*}\mu^{[\gamma_{v}]}\end{equation}
 
  Taking   smooth parts of our maps above gives the following diagram.
\[\begin{tikzcd}\Msw_{g,n}(B)\ar{rr} & & B^{n}\times \bar M_{g-1,n+2}/\mathbb Z_{2}\subset B^{n}\times \bar M_{g,n}
\\ & B^{n+1}\times \bar M_{g-1,n+2}\dar{\totl{\Delta}} \rar{\totl{\pi}}& B^{n}\times \bar M_{g-1,n+2}\uar{\totl{\psi}}
\\ \Msw_{g-1,n+2}(B)\rar & B^{n+2}\times \bar M_{g-1,n+2}\end{tikzcd}\]
Let $C_{g-1,n+2}$ be the pushforward of $\mu^{[\gamma_{v}]}$ in $H^{*}(B^{n+2}\times \bar M_{g-1,n+2})$ using the smooth part map, $\mathcal M_{g-1,n+2}\longrightarrow \bar M_{g-1,n+2}$, and let  $C_{g,n}$ in $H^{*}(B^{n}\times\bar M_{g,n})$ be the pushforward of $\mu$.   Konsevich and Mannin's genus-reduction formula from \cite{KM} can be restated\footnote{Actually, the  genus-reduction axiom from \cite{KM} also keeps track of the homology class of curves, whereas our formula only keeps track of their $\omega$--energy. In this case, our formula can easily be upgraded to keep track of homology classes as outlined in section \ref{fgf}.} as \[\totl{\psi}^{*} C_{g,n}=\totl{\pi}_{!}\totl{\Delta}^{*}C_{g-1,n+2}\]
which is implied by our gluing formula, (\ref{kmgf}), because each stage of pushing forward or pulling back commutes with pushing forward using the smooth part map.

The splitting axiom is proved similarly, except  now $\gamma$ is a tropical curve with $2$ vertices connected to one internal edge, and $n_{i}$ external edges attached to the $i$th vertex. Now $\gamma$  has no symmetries (assuming $n\neq 0$).  The relevant commutative diagram of maps for genus $g$ invariants is below. (The disjoint unions below are over choices of nonnegative integers $g_{i}$ so that $g_{1}+g_{2}=g$.) 
\[\begin{tikzcd}B^{n}\times\mathcal M_{g,n}\rar{\totl\cdot} &  B^{n}\times \bar M_{g,n}
\\B^{n}\times \mathcal M_{g,n}\tc\gamma\rar{\totl\cdot}&B^{n}\times \coprod\bar M_{g_{1},n_{1}+1}\times \bar M_{g_{2},n_{2}+1}\uar[hook]{\psi}
\\  \dar{\Delta}B^{n+1}\times \ex T\rtimes \coprod\mathcal M_{g_{1},n_{1}+1}\times \mathcal M_{g_{2},n_{2}+1}\uar{I^{[\gamma]}}\rar{\totl\cdot} &
B^{n+1}\times \coprod\bar M_{g_{1},n_{1}+1}\times \bar M_{g_{2},n_{2}+1}\dar{\totl\Delta}\uar{\pi}
\\ B^{n+2}\times \coprod\mathcal M_{g_{1},n_{1}+1}\times \mathcal M_{g_{2},n_{2}+1}\rar{\totl\cdot}& B^{n+2}\times \coprod\bar M_{g_{1},n_{1}+1}\times \bar M_{g_{2},n_{2}+1} 
\end{tikzcd}\]
Our gluing formula, $\mu\tc \gamma=I^{[\gamma]}_{!}\Delta^{*}(\mu^{[\gamma_{v_{1}}]}\wedge\mu^{[\gamma_{v_{2}}]})$, implies 
\[\psi^{*}C_{g,n}=\sum_{g_{1}+g_{2}=g}\pi_{!}\totl{\Delta}^{*}(C_{g_{1},n_{1}}\wedge C_{g_{2},n_{2}})\]
which is the splitting axiom from \cite{KM}.

\section{Further gluing formulae}\label{fgf}

We now have two gluing formulae. The first is
\[\eta\tc \gamma=\frac {k_{\gamma}}{\abs{\Aut \gamma}}i^{[\gamma]}_{!}\Delta^{*}\prod_{v}\eta^{[\gamma_{v}]}\]
 where $k_{\gamma}$ is the product of the multiplicities of the internal edges of $\gamma$, and $\eta$ and $\eta^{[\gamma_{v}]}$ are the pushforward of $q^{E}\hbar^{2g-2+n}$ via the maps
 \[ev\co  \Msw_{\cdot}\ex B\longrightarrow\coprod_{n} (\rend\ex B)^{n}\]
 \[ev^{[\gamma_{v}]}\co  \Msw_{[\gamma_{v}]}(\ex B\tc v)\longrightarrow \rend^{[\gamma_{v}]}\ex B\]
 and $\Delta$ and $i^{[\gamma]}$ are natural maps
 \[\begin{tikzcd}\prod_{v}\rend^{[\gamma_{v}]}(\ex B\tc v)&\lar{\Delta}\rend^{[\gamma]}(\ex B)\tc{\totb{ev^{[\gamma]}}\gamma}\rar{i^{[\gamma]}}&(\coprod_{n}(\rend\ex B)^{n})\tc{\totb{ev}\gamma}\end{tikzcd}\ .\]
 Our second gluing formula
 \[\mu\tc \gamma=\frac {1}{\abs{\Aut \gamma}}I^{[\gamma]}_{!}\Delta^{*}\prod_{v}\mu^{[\gamma_{v}]}\]
involves the pushforward, $\mu$ and $\mu^{[\gamma_{v}]}$, of $q^{E}$ using enhanced evaluation maps
 \[EV\co  \Msw_{\cdot}\ex B\longrightarrow \mathcal X\ex B\]
 \[EV^{[\gamma_{v}]}\co  \Msw_{[\gamma_{v}]}(\ex B\tc v)\longrightarrow \mathcal X^{[\gamma_{v}]}(\ex B\tc v)\]
 and the natural maps
 \[\begin{tikzcd}\prod_{v}\mathcal X^{[\gamma_{v}]}(\ex B\tc v)&\lar{\Delta}(\mathcal X^{[\gamma]}\ex B)\tc{\totb{EV^{[\gamma]}}\gamma}\rar{I^{[\gamma]}}&(\mathcal X\ex B)\tc{\totb{EV}\gamma}\end{tikzcd}\ .\]

Our first gluing formula follows from Theorem \ref{sg} and Lemma \ref{eta pf}, and the second from Theorem \ref{xg} and Lemma \ref{xpf}.  These two key theorems follow from the fact that the left two squares in the following diagram become fiber-product diagrams after applying tropical completion suitably, as in Lemmas \ref{tcfp} and \ref{xpf}.
\[\begin{tikzcd}\Msw_{[\gamma]}(\ex B)\dar{cut}\rar{EV^{[\gamma]}}\ar[bend left, swap]{rr}{\ev^{[\gamma]}}\ar[bend left]{rrr}{ev^{[\gamma]}}&\mathcal X^{[\gamma]}\rar\dar{\Delta}&\End^{[\gamma]}\ex B\rar\dar{\Delta}&\rend^{[\gamma]}\ex B\dar{\Delta}
\\ \prod_{v}\Msw_{[\gamma_{v}]}(\ex B\tc v)\rar{\prod_{v}EV^{[\gamma_{v}]}}\ar[bend right]{rr}{\prod_{v}\ev^{[\gamma_{v}]}}\ar[bend right,swap]{rrr}{\prod_{v}ev^{[\gamma_{v}]}}&\prod_{v}\mathcal X^{[\gamma_{v}]}(\ex B\tc v)\rar&\prod_{v}\End^{[\gamma_{v}]}(\ex B\tc v)\rar& \prod_{v}\rend^{[\gamma_{v}]}(\ex B\tc v)\end{tikzcd}\] 
We can also construct gluing formulae keeping track of more discrete information. Let $\hat{\mathcal X}$ and $\hat{\mathcal X}^{[\gamma_{v}]}$ be  covers of $\mathcal X$ and $\mathcal X^{[\gamma_{v}]}$ respectively  with lifts $\hat {EV}$ and $\hat{EV}^{[\gamma]}$ of our evaluation maps.
\[\begin{tikzcd}& \hat{\mathcal X}\dar & & \hat{\mathcal X}^{[\gamma_{v}]}\dar 
\\ \Msw(\ex B)\rar{EV}\ar{ur}{\hat{EV}}& \mathcal X &\Msw_{[\gamma_{v}]}(\ex B\tc v)\rar{EV}\ar{ur}{\hat{EV}^{[\gamma_{v}]}}& \mathcal X^{[\gamma_{v}]}\end{tikzcd}\]
Suppose further that there is a lift, $\hat I^{[\gamma]}$ of $I^{[\gamma]}$ as in the diagram below, compatible with $\hat{EV}$ and $\hat{EV}^{[\gamma_{v}]}$.
\[\begin{tikzcd} \prod_{v}\hat{\mathcal X}^{[\gamma_{v}]}\dar&\lar{\Delta}(\prod_{v}\hat {\mathcal X}^{[\gamma_{v}]})\times_{\prod_{v}\mathcal X^{[\gamma_{v}]}}\mathcal X^{[\gamma]}\tc{\totb{EV}\gamma} \dar \rar{\hat I^{[\gamma]}} & \hat{\mathcal X}\tc{\totb{\hat{EV}}\gamma}\dar
\\\prod_{v}\mathcal X^{[\gamma_{v}]} &\lar{\Delta}  \mathcal X^{[\gamma]}\tc{\totb{EV}\gamma}\rar{I^{[\gamma]}}&\mathcal X\tc{\totb{EV}\gamma}
\end{tikzcd}\]
The following gluing formula then holds,  
\[\hat\mu\tc \gamma=\frac 1{\abs{\Aut \gamma}}\hat I^{[\gamma]}_{!}\Delta^{*}\prod_{v}\hat \mu^{[\gamma_{v}]} \]
where 
\[\hat\mu:=\hat{EV}_{!}(q^{E}) \text{ and }\hat \mu^{[\gamma_{v}]}:=\hat{EV}^{[\gamma_{v}]}_{!}(q^{E})\ .\]
We can also include the contribution of gravitational descendants, as in Theorem \ref{xg} and Lemma \ref{xpf}.

An example of such a lift is as follows:  we can lift $EV$ to keep track of the integral over curves  of all closed $2$--forms $\alpha$ in $\Omega^{2}(\ex B)$, and lift $EV^{[\gamma_{v}]}$ to record the integral of $\alpha\tc v$ over curves.

Of course, we can similarly extend our first gluing formula by suitably lifting $ev$. An example of this construction is the use of rim tori in \cite{IP}, explained further in \cite{ZingerRT}.

 \bibliographystyle{plain}

\bibliography{ref.bib}
 \end{document}